\documentclass[twoside,11pt]{article}

\usepackage{jmlr2e}

  \usepackage{hyperref,url,breakurl,fancybox,amsmath,bm,paralist,booktabs,subfig,graphicx}%fancybox => rounded box, breakurl => url-s can break
  \usepackage[mathscr]{eucal}
  \usepackage{enumitem,color} %customizing list: leftmargin, ...
  \newcommand{\X}{\mathscr{X}}   %domain of the distributions
  \renewcommand{\H}{\mathscr{H}} %RKHS(K)
  \newcommand{\Hmeta}{\mathscr{G}} %RKHS-meta
  \newcommand{\Kmeta}{K_{\Hmeta}} %K-meta
  \newcommand{\R}{\mathbb{R}}    %real number
  \newcommand{\Z}{\mathbb{Z}}%   integer
  \newcommand{\E}{\mathbb{E}}    %expectation
  \renewcommand{\Pr}{\mathbb{P}} %probability
  \renewcommand{\L}{\mathscr{L}} %linear operator
  \newcommand{\Eo}{\mathcal{R}}  %expected risk
  \newcommand{\EoS}{\mathcal{E}}  %excess risk
  \renewcommand{\d}{\mathrm{d}}  %d(x): integration
  \renewcommand{\b}{\mathbf}     %bold
  \newcommand{\A}{\mathscr{A}}%CDV-1
  \newcommand{\B}{\mathscr{B}}%CDV-2
  \newcommand{\N}{\mathscr{N}}%CDV-3
  \newcommand{\M}{\mathscr{M}}%probability measure
  \renewcommand{\P}{\mathcal{P}}%P(b,c)
  \renewcommand{\O}{\mathcal{O}}%Ordo
  \newcommand{\tb}{\textbf}%bold text
  \newcommand{\Bo}{\mathcal{B}}%Borel sets
  \newcommand{\I}{\mathbb{I}}%indicator function
  \renewcommand{\S}{\mathscr{S}}%sigma-alg

  \newcommand{\cl}[1]{cl\big[#1\big]}%closure of a subspace
  \renewcommand{\Im}{Im}%image of an operator
  \newcommand{\Ker}{Ker}%nullspace of an operator
  \newcommand{\K}{\mathscr{K}}%general Hilbert space
  \newcommand{\G}{\mathscr{G}}

  \DeclareMathOperator*{\argmin}{arg\,min}

\jmlrheading{17}{2016}{1-40}{12/14; Revised 1/16}{9/16}{Zolt{\'a}n Szab{\'o}, Bharath K. Sriperumbudur, Barnab{\'a}s P{\'o}czos, and Arthur Gretton}

\ShortHeadings{Learning Theory for Distribution Regression}{Szab{\'o} et al.}%<=50chars
\firstpageno{1}

\begin{document}

\title{Learning Theory for Distribution Regression}
\author{\name Zolt{\'a}n Szab{\'o}\thanks{Now at Department of Applied Mathematics, CMAP, {\'E}cole Polytechnique, Route de Saclay, 91128 Palaiseau, France.} 
        \email zoltan.szabo@gatsby.ucl.ac.uk\\
        ORCID 0000-0001-6183-7603\\
        \addr
        Gatsby Unit, University College London\\
        Sainsbury Wellcome Centre, 25 Howland Street\\
	London - W1T 4JG, UK
	\AND
	\name Bharath K. Sriperumbudur \email bks18@psu.edu \\
        \addr Department of Statistics\\
        Pennsylvania State University\\
        University Park, PA 16802, USA
        \AND
        \name Barnab{\'a}s P{\'o}czos \email bapoczos@cs.cmu.edu \\
        \addr Machine Learning Department\\
        School of Computer Science\\
        Carnegie Mellon University\\
	5000 Forbes Avenue Pittsburgh PA 15213 USA
	\AND
	\name Arthur Gretton
        \email arthur.gretton@gmail.com\\
        	ORCID 0000-0003-3169-7624\\
        \addr
        Gatsby Unit, University College London\\
        Sainsbury Wellcome Centre, 25 Howland Street\\
	London - W1T 4JG, UK
}
\editor{Ingo Steinwart}
\maketitle

\begin{abstract}%   <- trailing '%' for backward compatibility of .sty file; <=200 words
We focus on the distribution regression problem: regressing to vector-valued outputs from probability measures.
Many important machine learning and statistical tasks fit into this framework, including multi-instance learning and point estimation problems without analytical solution (such as hyperparameter or entropy estimation). 
Despite the large number of available heuristics in the literature, the inherent two-stage sampled nature of the problem 
makes the theoretical analysis quite challenging, since 
in practice only samples from sampled distributions are observable, and the estimates have to rely on  similarities computed between sets of points. To the best of our knowledge, the only existing technique with 
consistency guarantees for distribution regression requires kernel density estimation as an intermediate step (which often performs poorly in practice), and the domain of the distributions to be
compact Euclidean. In this paper, we study a simple, analytically computable, ridge regression-based alternative to distribution regression, where we embed the distributions to a reproducing kernel Hilbert space, and learn
the regressor from the embeddings to the outputs. Our main contribution is to prove that this scheme is consistent in the two-stage sampled setup under mild conditions (on separable topological domains enriched with kernels):
we present an exact computational-statistical efficiency trade-off analysis showing that our estimator is able to match the \emph{one-stage} sampled minimax optimal rate \citep{caponnetto07optimal,steinwart09optimal}. This result answers a $17$-year-old open question, establishing 
the consistency of the classical set kernel \citep{haussler99convolution,gartner02multi} in regression.  We also  cover consistency for more recent kernels on distributions, 
including those due to \citet{christmann10universal}.
\end{abstract}


\section{Introduction}\label{sec:intro}
We  address the learning problem of \emph{distribution regression} in the two-stage sampled setting, where we only have bags of samples from the probability distributions: 
we  regress from probability measures to real-valued  \citep{poczos13distribution} responses, or more generally to vector-valued outputs (belonging to an arbitrary separable Hilbert space).
Many classical problems in machine learning and statistics can be analysed in this framework. On the machine learning side, multiple instance learning \citep{dietterich97solving,ray01multiple,dooly02multiple} can be thought of in this way, 
where each instance in a labeled bag is an i.i.d.\ (independent identically distributed) sample from a distribution. On the statistical side, tasks might include point estimation of statistics on a 
distribution without closed form analytical expressions (e.g., its entropy or a hyperparameter).

\tb{Intuitive description of our goal:} Let us start with a somewhat informal definition of the distribution regression problem and an intuitive phrasing of 
our goals. Suppose that our data consist of $\b{z}=\{(x_i,y_i)\}_{i=1}^l$, where $x_i$ is a probability distribution, $y_i$ is its label (in the simplest case $y_i\in \R$ or $y_i\in\R^d$) and
each $(x_i,y_i)$  pair is i.i.d.\ sampled from a meta distribution $\mathcal{M}$. However, we do not observe $x_i$ directly; rather, we observe a sample
$x_{i,1},\ldots, x_{i,N_i} \stackrel{i.i.d.} {\sim} x_i$. 
Thus the observed data are $\hat{\b{z}} = \{(\{x_{i,n}\}_{n=1}^{N_i},y_i)\}_{i=1}^l$. 
Since $x_{i,j}$ is sampled from $x_i$, and $x_i$ 
is sampled from $\mathcal{M}$, we call this process two-stage sampling.
Our goal is to predict a new $y_{l+1}$ from a new batch of samples $x_{l+1,1},\ldots, x_{l+1,N_{l+1}}$ drawn from a new distribution $x_{l+1}\sim \mathcal{M}$. 
For example, in a medical application, the $i^{th}$ patient might be identified with a probability distribution ($x_i$), which can be periodically assessed by
blood tests ($\{x_{i,n}\}_{n=1}^{N_i}$). We are also given some health indicator of the patient ($y_i$), which might be inferred from his/her blood measurements.
Based on the observations ($\hat{\b{z}}$), we might wish to learn the mapping from the set of blood tests to the health indicator, and the hope is that by 
observing more patients (larger $l$) and performing a larger number of tests (larger $N_i$) the estimated mapping ($\hat{f}=\hat{f}({\hat{\b{z}}})$) becomes more ``precise''. Briefly, we consider the following questions:
\begin{center}
\fbox{\begin{minipage}{27.8em}
  Can the distribution regression problem be solved consistently under mild conditions? What is the exact computational-statistical efficiency trade-off implied by the two-stage sampling?
\end{minipage}}
\end{center}
In our work the estimated mapping ($\hat{f}$) is the analytical solution of a kernel ridge regression (KRR) problem.\footnote{Beyond its simple analytical formula, kernel ridge regression also allows efficient distributed \citep{zhang15divide,richtarik16distributed}, sketch \citep{alaoui15fast,yang16randomized} and Nystr{\"o}m based approximations \citep{rudi15less}.}
 The performance of $\hat{f}$ depends on the assumed function class ($\H$), 
the family of $\hat{f}$ candidates used in the ridge formulation. We shall focus on the analysis of two settings:
\begin{enumerate}
 \item \tb{Well-specified case ($f_*\in\H$):}\label{eq:f*-def} In this case we assume that the regression function $f_*$ belongs to $\H$. 
	We focus on bounding the goodness of $\hat{f}$ compared to $f_*$. In other words, if $\Eo[f_*]$ denotes the prediction error (expected risk) of $f_*$, then our goal is to
	derive a finite-sample bound for the excess risk, $\EoS(\hat{f},f_*)=\Eo[\hat{f}] - \Eo[f_*]$ that holds with high probability.
	We make use of this bound to establish the consistency of the estimator (i.e., drive the excess risk to zero) and to derive the exact computational-statistical efficiency trade-off of the estimator 
        as a function of the sample number ($l$, $N=N_i,\,\forall\, i$) and the problem difficulty (see Theorem~\ref{conseq:conv-rate} and its corresponding remarks for more details).
  \item \tb{Misspecified case ($f_*\in L^2\backslash \H$):}\label{motiv:misspec} Since in practise it might be hard to check whether $f_*\in \H$, we also study the misspecified setting of $f_*\in L^2$; 
	  the relevant case is when $f_*\in L^2\backslash \H$. In the misspecified setting the 'richness' of $\H$ has crucial importance, in other words the size of $D_{\H}^2 = \inf_{f\in \H}\|f_*-f\|_2^2$, the approximation error from $\H$. 
          Our main contributions consist of proving a finite-sample excess risk bound, using which 
          we show that the proposed estimator can attain the ideal performance, i.e., $\EoS(\hat{f},f_*) - D_{\H}^2$ can be driven to zero. Moreover, 
	  on smooth classes of $f_*$-s, we give a simple and explicit description for the computational-statistical efficiency trade-off of our estimator (see Theorem~\ref{conseq:L2rate:Assumption2b} and its corresponding remarks for more details).
\end{enumerate}

There exist a vast number of heuristics to tackle learning problems on distributions; we will review them in Section~\ref{sec:related-work}. However, to the best of our knowledge, the only prior work 
addressing the \emph{consistency} of regression on distributions requires kernel density estimation
\citep{poczos13distribution,oliva14fast,sutherland16lineartime}, which assumes that the response variable is 
scalar-valued,\footnote{\citet{oliva13ICML,oliva15fast} consider the case where the responses are also distributions or functions.} and the covariates are nonparametric continuous distributions on $\R^d$.
 As in our setting, the exact forms of these distributions are unknown; they are available only through finite sample sets. P{\'o}czos et al.\ estimated these distributions
through a kernel density estimator (assuming these distributions have a density) and then constructed a kernel regressor that acts on these kernel
density estimates.\footnote{We would like to clarify that the kernels used in their work are classical smoothing kernels---extensively studied in non-parametric statistics \citep{gyorfi2002}---and not the reproducing kernels that appear throughout our paper.} 
Using the classical bias-variance decomposition analysis for kernel regressors, they showed the consistency of the constructed kernel regressor, and provided a polynomial upper bound 
on the rates, assuming the true regressor to be H{\"o}lder continuous, and the meta distribution that generates the covariates $x_i$ to have finite doubling 
dimension \citep{kpotufe2011}.\footnote{Using a random kitchen sinks approach, with orthonormal basis 
projection estimators \citet{oliva14fast,sutherland16lineartime} propose distribution regression algorithms 
that can computationally handle large scale datasets; 
as with \citet{poczos13distribution}, these approaches are based on density estimation in $\R^d$.\label{footnote:RNDkitchensink}}
 
One can define kernel learning algorithms on bags based on set kernels \citep{gartner02multi} by computing the similarity of the sets/bags of samples representing the input distributions; set kernels are also called called multi-instance kernels or ensemble kernels, and are examples of convolution kernels \citep{haussler99convolution}. 
In this case, the similarity of two sets is measured by the average pairwise point similarities between the sets. From a theoretical perspective,
nothing is known about the consistency of set kernel based learning method since their introduction in 1999 \citep{haussler99convolution,gartner02multi}:
i.e.\ in what sense (and with what rates) is the learning algorithm consistent, when the number of items per bag, and the number of bags, are allowed to increase?

It is possible, however, to view set kernels in a distribution setting, as they represent valid kernels between (mean) embeddings of empirical probability measures into a reproducing kernel Hilbert space \citep[RKHS;][]{BerTho04}.
The population limits are well-defined as being dot products between the embeddings of the generating distributions \citep{altun06unifying}, and for characteristic kernels the distance 
between embeddings defines a metric on probability measures \citep{sriperumbudur11universality,gretton12kernel}.  When bounded kernels are used, mean embeddings exist for all 
probability measures \citep{fukumizu04dimensionality}. When we consider the distribution regression setting, however, there is no reason to limit ourselves to set kernels.
Embeddings of probability measures to RKHS are  used by \citet{christmann10universal} in defining a yet larger class
of easily computable kernels on distributions, via operations performed on the embeddings and their distances.
Note that the relation between set kernels and kernels on distributions was also applied by \citet{muandet12learning} for classification on distribution-valued inputs, however consistency was not studied in that work.
We also note that motivated by the current paper, \citet{lopez15towards} have recently presented the first theoretical results about surrogate risk guarantees on a class 
(relying on uniformly bounded Lipschitz functionals) of soft distribution-classification problems.

Our \tb{contribution} in this paper is to establish the learning theory of a simple, mean embedding based ridge regression (MERR) method for the distribution regression problem. 
This result applies both to the basic set kernels of \citet{haussler99convolution,gartner02multi}, the distribution kernels of  \citet{christmann10universal}, and additional related kernels.
We provide finite-sample excess risk bounds, prove consistency, and show how the two-stage sampled nature of the problem (bag size) governs the computational-statistical efficiency of the MERR estimator.
More specifically, in the \begin{enumerate}[labelindent=0cm,leftmargin=*,topsep=0cm,partopsep=0cm,parsep=0cm,itemsep=0cm]
    \item \tb{well-specified case:} We\vspace{1mm}
	  \begin{enumerate}[labelindent=0cm,leftmargin=*,topsep=0cm,partopsep=0cm,parsep=0cm,itemsep=0cm]
	      \item derive finite-sample bounds on the excess risk: We construct $\Eo[\hat{f}] - \Eo[f_*]\le r(l,N,\lambda)$ bounds holding with high probability, 
		    where $\lambda$ is the regularization parameter in the ridge problem ($\lambda\rightarrow 0$, $l\rightarrow \infty$, $N=N_i\rightarrow \infty$). \vspace{1mm}
	      \item establish consistency and computational-statistical efficiency trade-off of the MERR estimator
		      on a general prior family $\P(b,c)$ as defined by \citet{caponnetto07optimal}, where $b$ captures the effective input dimension, and larger $c$  means 
		    smoother $f_*$ ($1<b$, $c\in (1,2]$).
 In particular, when the number of samples per bag is chosen as $ N=l^{a}\log(l)$ and  $a\ge \frac{b(c+1)}{bc+1}$, then the learning rate saturates at $l^{-\frac{bc}{bc+1}}$, which is known to be one-stage sampled minimax optimal \citep{caponnetto07optimal}.
In other words, by choosing  $a =  \frac{b(c+1)}{bc+1} < 2$, we suffer
{\em no loss in statistical performance} compared with the {\em best possible one-stage sampled estimator}.\vspace{1mm}  
	            \end{enumerate}
	  Note: the advantage of considering the  $\P(b,c)$ family is two-fold. It does not assume parametric distributions, yet certain complexity terms can be explicitly upper bounded in the family. This property will be exploited in our analysis.
		Moreover, (for special input distributions) the parameter $b$ can be related to the spectral decay of Gaussian Gram matrices, and existing analysis techniques \citep{steinwart08support} 
		may be used in interpreting these decay conditions.
    \item \tb{misspecified case:}   We establish consistency and convergence rates even if $f_*\notin \H$. Particularly, by deriving finite-sample bounds on the excess risk we\vspace{1mm}
	  \begin{enumerate}[labelindent=0cm,leftmargin=*,topsep=0cm,partopsep=0cm,parsep=0cm,itemsep=0cm]
		\item prove that the MERR estimator can achieve the best possible approximation accuracy from $\H$, i.e. the $\Eo[\hat{f}] - \Eo[f_*] - D_{\H}^2$ quantity can be driven to zero (recall that $D_{\H} = \inf_{f\in \H}\|f_*-f\|_2$). Specifically, this result implies that if $\H$ is dense in $L^2$ ($D_{\H}=0$), then the excess risk $\Eo[\hat{f}] - \Eo[f_*]$ converges to zero.\vspace{1mm}
	        \item analyse the computational-statistical efficiency trade-off: We show that by choosing the bag size 
		  as $N=l^{2a}\log(l)$ $(a>0)$ one can get rate $l^{-\frac{2sa}{s+1}}$ for $a\le \frac{s+1}{s+2}$, and the rate saturates for $a\ge \frac{s+1}{s+2}$ at 
		  $l^{-\frac{2s}{s+2}}$, where the difficulty of the regression problem is captured by $s\in(0,1]$ (a larger $s$ means an easier problem). 
		  This means that  easier tasks give rise to faster convergence (for $s=1$, the rate is $l^{-\frac{2}{3}}$),  the bag size $N$ can again be \emph{sub-quadratic} in $l$ ($2a\le\frac{2(s+1)}{s+2}\le \frac{4}{3}<2$),  and 
		  the rate at saturation is close to  $\tilde{r}(l)=l^{-\frac{2s}{2s+1}}$, which is the asymptotically optimal rate in the one-stage sampled setup, with real-valued output 
		  and stricter eigenvalue decay conditions \citep{steinwart09optimal}.
	  \end{enumerate}
\end{enumerate}

Due to the differences in the assumptions made and the loss function used, a direct comparison of our theoretical result and that of \citet{poczos13distribution} 
remains an open question, 
however we make three observations. First, our approach is more general, since we may regress from any probability measure defined on separable, topological domains endowed with kernels.
P{\'o}czos et al.'s work is restricted to compact domains of finite dimensional Euclidean spaces, and requires the distributions to admit probability densities; distributions on strings, 
graphs, and other structured objects are disallowed. Second, in our analysis we will allow separable Hilbert space valued outputs, in contrast to the real-valued output considered by \citet{poczos13distribution}. Third, density estimates in high dimensional spaces suffer from slow convergence rates \citep[Section 6.5]{Wasserman06}. Our approach 
mitigates this problem, as it works directly on distribution embeddings, and does not make use of density estimation as an intermediate step. 

The principal challenge in proving theoretical guarantees arises from the two-stage sampled nature of the inputs. In our analysis of the well-specified case, we make use of \citet{caponnetto07optimal}'s results, 
which focus (only) on the one-stage sample setup. These results will make our analysis somewhat shorter (but still rather challenging) by giving upper bounds for some of the  objective terms. 
Even the verification of these conditions requires care since the inputs in the ridge regression are themselves distribution embeddings (i.e., functions in a reproducing kernel Hilbert space).

\label{paragraph:misspec} In the misspecified case,  RKHS methods alone are not sufficient to obtain excess risk bounds: one has to take into account the ``richness'' of the modelling RKHS class ($\H$) in the embedding $L^2$ space.
The fundamental challenge  is whether it is possible to achieve the best possible performance dictated by $\H$; or in the special case when further smoothness conditions hold on 
$f_{*}$, what convergence rates can yet be attained, and what computational-statistical efficiency trade-off realized. The second smoothness property could be modelled for example by range spaces of (fractional) powers of integral operators associated to $\H$. Indeed, 
there exist several results along these lines with KRR for the case of real-valued outputs: see for example 
\citep[Theorem~1.1]{sun09application}, \citep[Corollary~3.2]{sun09note}, \citep[Theorem 3.7 with Assumption 3.2]{mendelson10regularization}. The question of optimal rates has also been addressed for
the semi-supervised KRR setting \citep[Theorem~1]{caponnetto06optimal}, and for clipped KRR estimators  \citep{steinwart09optimal} with integral operators of rapidly decaying spectrum.
Our results apply more generally to the two-stage sampled setting and to vector valued outputs belonging to separable Hilbert spaces. Moreover,
we obtain a general consistency result without range space assumptions, showing that the modelling power of $\H$ can be fully exploited, and convergence to 
the best approximation available from $\H$ can be realized.\footnote{Specializing our result, we get explicit rates and an exact computational-statistical efficiency description for MERR as a function of sample numbers and problem difficulty, for smooth regression functions.}

There are numerous areas in machine learning and statistics, where estimating vector-valued functions has crucial importance. Often in statistics,  one is not only confronted with the estimation of a scalar parameter, 
but with a vector of parameters.  On the machine learning side, multi-task  learning \citep{evgeniou05learning}, functional response regression \citep{kadri16operator}, or structured output prediction \citep{brouard11semi,kadri13generalized} fall under the same umbrella: they can be naturally phrased as
learning vector-valued functions \citep{micchelli05learning}.  The idea  underlying all these tasks is simple and intuitive: if multiple prediction problems have to be solved simultaneously, it might be beneficial to
exploit their dependencies. Imagine for example that the task is to predict the motion of a dancer: taking into account the interrelation of the actor's body parts is likely to lead to more accurate estimation, as opposed to predicting the individual parts one by one, independently.
Successful real-world applications of a multi-task approach include for example preference modelling of users with similar demographics \citep{evgeniou05learning}, prediction of the daily precipitation profiles of weather stations \citep{kadri10nonlinear}, 
acoustic-to-articulatory speech inversion \citep{kadri16operator}, identifying biomarkers capable of 
tracking the progress of Alzheimer's disease \citep{zhou13modeling}, 
personalized human activity recognition based on iPod/iPhone accelerometer data \citep{sun13large}, finger trajectory prediction in brain-computer interfaces \citep{kadri12multiple} or ecological 
inference \citep{flaxman15who}; for a recent review on multi-output prediction methods see \citep{alvarez11kernels,borchani15survey}. A mathematically sound way of encoding prior information about the relation of the outputs can be realized by operator-valued kernels and the associated vector-valued RKHS-s \citep{pedrick57theory,micchelli05learning,carmeli06vector,carmeli10vector}; 
this is the tool we use to allow vector-valued learning tasks.

Finally, we note that the current work extends  our earlier conference paper \citep{szabo15twostage} in several important respects: 
we now show that the MERR method can attain the one-stage sampled minimax optimal rate; 
we  generalize the analysis in the well-specified setting to allow outputs belonging to an arbitrary separable Hilbert spaces (in contrast to the original scalar-valued output domain); and we tackle the misspecified setting, obtaining finite sample guarantees, consistency, and computational-statistical efficiency trade-offs.

The paper is structured as follows: The distribution regression problem and the MERR technique are introduced in Section \ref{sec:problem}.
Our assumptions are detailed in Section \ref{sec:assumptions}. We present our theoretical guarantees (finite-sample bounds on the excess risk, consistency, computational-statistical efficiency trade-offs) 
in Section~\ref{sec:convergence analysis}: the well-specified case is considered in Section~\ref{sec:results:A1}, and the misspecified setting is the focus of Section~\ref{sec:results:A2}. 
Section~\ref{sec:related-work} is devoted to an overview of existing heuristics for learning on distributions.
Conclusions are drawn in Section~\ref{sec:conclusions}. Section~\ref{sec:proofs} contains proof details. In Section~\ref{sec:assumptions-discussion} we discuss our assumptions with concrete examples. 

\section{The Distribution Regression Problem} \label{sec:problem}
Below we first introduce our notation (Section~\ref{sec:notations}), then formally define the distribution regression task (Section~\ref{sec:DR problem}).

\subsection{Notation}\label{sec:notations} We use the following notations throughout the paper:
\begin{itemize}[labelindent=0cm,leftmargin=*,topsep=0cm,partopsep=0cm,parsep=0cm,itemsep=0cm]
 \item 
      \tb{Sets, topology, measure theory:} 
      Let $\K$ be a Hilbert space; $\cl{V}$ is the closure of a set $V\subseteq \K$. $\times_{i\in I} S_i$ is the direct product of sets $S_i$.
      $f\circ g$ is the composition of function $f$ and $g$. Let $(\X,\tau)$ be a topological space and let $\Bo(\X):=\Bo(\tau)$ be the Borel $\sigma$-algebra induced by the topology $\tau$. 
      If $(\X,d)$ is a metric space, then $\Bo=\Bo(d)$ is the Borel $\sigma$-algebra generated by the open sets induced by metric $d$. $\M_1^+(\X)$ denotes the set of Borel probability measures on the $(\X,\Bo(\X))$ measurable space. 
      Given measurable spaces $(U_1,\S_1)$ and $(U_2,\S_2)$,  the $\S_1\otimes \S_2$ product $\sigma$-algebra \citep[page~480]{steinwart08support} on the product space $U_1 \times U_2$ is the $\sigma$-algebra generated by the cylinder sets $U_1\times S_2$, $S_1\times U_2$ ($S_1\in \S_1$, $S_2\in \S_2$). 
      The weak topology ($\tau_w=\tau_w(\X,\tau)$) on $\M^+_1(\X)$ is defined as the weakest topology such that the $L_h:(\M^+_1(\X),\tau_w) \rightarrow \R$, 
      $L_h(x)=\int_{\X}h(u)\d x(u)$ mapping is continuous for all $h\in C_b(\X)=\{(\X,\tau)\rightarrow  \R\text{ bounded, continuous functions}\}$.\vspace{1mm}
 \item
      \tb{Functional analysis:} Let $(N_1,\left\|\cdot\right\|_{N_1})$ and $(N_2,\left\|\cdot\right\|_{N_2})$ denote two normed spaces, then  $\L(N_1,N_2)$ stands for the 
      space of $N_1\rightarrow N_2$ bounded linear operators; if $N_1=N_2$, we will use the $\L(N_1)=\L(N_1,N_2)$ shorthand. For $M\in\L(N_1,N_2)$ the operator norm is defined as 
      $\left\|M\right\|_{\L(N_1,N_2)}=\sup_{0\ne h\in N_1}\left\|Mh\right\|_{N_2}/\left\|h\right\|_{N_1}$, 
       $\Im(M)=\{Mn_1\}_{n_1\in N_1}$ denotes the range of $M$, $\Ker(M)=\{n_1\in N_1: Mn_1 = 0\}$ is 
      the null space of $M$. 
      Let $\K$ be a Hilbert space. 
      The adjoint operator $M^*\in \L(\K)$ of an operator $M \in \L(\K)$ is the operator such that
      $\left<Ma,b\right>_{\K} = \left<a,M^*b\right>_{\K}$ for all $a$ and $b$ in $\K$.      
      $M\in \L(\K)$ is called  positive if $\left<Ma,a\right>_{\K}\ge 0$ ($\forall a\in\K$),  self-adjoint if $M=M^*$, and trace class  if 
      $\sum_{j\in J}\left<|M|e_j,e_j\right>_{\K}<\infty$ for an $(e_j)_{j\in J}$ ONB (orthonormal basis) of $\K$ ($|M|:=(M^*M)^{\frac{1}{2}}$), in which case  $Tr(M) := \sum_{j\in J}\left<Me_j,e_j\right>_{\K}<\infty$; 
       compact if $cl\left[Ma: a\in\K, \left\|a\right\|_{\K}\le 1\right]$ is a compact set. 
      Let $\K_1$ and $\K_2$ be Hilbert spaces. $M\in \L(\K_1,\K_2)$ is called Hilbert-Schmidt if 
      $\left\|M\right\|_{\L_2(\K_1,\K_2)}^2 = Tr(M^*M) = \sum_{j\in J} \left<Me_j,Me_j\right>_{\K_2} < \infty$ 
      for some $(e_j)_{j\in J}$ ONB  of $\K_1$. The space of Hilbert-Schmidt operators is denoted by
      $\L_2(\K_1,\K_2)=\{M\in\L(\K_1,\K_2): \left\|M\right \|_{\L_2(\K_1,\K_2)}<\infty \}$. We use the shorthand notation $\L_2(\K)=\L_2(\K,\K)$ if $\K:=\K_1=\K_2$; $\L_2(\K)$ is separable if and only if $\K$ is separable \citep[page 506]{steinwart08support}. 
      Trace class and Hilbert-Schmidt operators over a $\K$ Hilbert space are compact operators \citep[page 505-506]{steinwart08support}; moreover,
      \begin{align}
	  \left\|A\right\|_{\L(\K)} &\le \left\|A\right\|_{\L_2(\K)}, \quad \forall A \in \L_2(\K),\label{eq:opnorm<=HSnorm}\\
	\left\|AB\right\|_{\L_2(\K)} &\le \left\|A\right\|_{\L_2(\K)} \left\| B \right\|_{\L(\K)}, \quad \forall A,B \in \L_2(\K).\label{eq:prod-norm-ineq}
      \end{align}
      $I$ is the identity operator; $I_l\in\R^{l\times l}$ is the identity matrix. \vspace{1mm}
\item 
    \tb{RKHS, mean embedding:} Let $H=H(k)$ be an RKHS \citep{steinwart08support} with $k:\X\times \X \rightarrow \R$ as the reproducing kernel. 
    Denote by 
    \begin{align*}
      X &= \mu\left(\M^+_1\left(\X\right)\right) = \{\mu_x: x\in \M^+_1\left(\X\right)\}\subseteq H, &   \mu_x&=\int_{\X}k(\cdot,u)\d x(u)=\E_{u\sim x}[k(\cdot,u)] \in H &
    \end{align*}
    the set of 
    mean embeddings \citep{BerTho04} of the distributions 
    to the space $H$.\footnote{The $x\mapsto \mu_x$ mapping is defined for \emph{all} $x\in \M^+_1(\X)$ if $k$ is bounded, i.e., $\sup_{u\in \X}k(u,u)<\infty$.\label{footnote:k-bounded-to-mu(X)}} 
    Let $Y$ be a separable Hilbert space, where the inner product is denoted by $\left<\cdot,\cdot\right>_Y$; the associated norm is $\left\|\cdot\right\|_Y$.
    $\H=\H(K)$ is the $Y$-valued RKHS \citep{pedrick57theory,micchelli05learning,carmeli06vector,carmeli10vector} of $X\rightarrow Y$ functions with $K:X\times X\rightarrow \L(Y)$
    as the reproducing kernel (we will present some concrete examples of $K$ in Section~\ref{sec:assumptions}; see Table~\ref{tab:K examples}); $K_{\mu_x}\in\L(Y,\H)$ is defined as
    \begin{align}
    K(\mu_x,\mu_t)(y) = (K_{\mu_t}y)(\mu_x),\quad (\forall \mu_x,\mu_t\in X),\text{ or } K(\cdot,\mu_t)(y) = K_{\mu_t}y \label{eq:K(a,b)-K_a}.
    \end{align}
    Further, $f(\mu_x) = K_{\mu_x}^*f$ $(\forall \mu_x \in X, f\in \H)$. \vspace{1mm}
\item 
    \tb{Regression function:} Let $\rho$ be the $\mu$-induced probability measure on the $Z=X\times Y$ product space, and let $\rho(\mu_x,y) = \rho(y|\mu_x)\rho_X(\mu_x)$ be the factorization of $\rho$ into 
    conditional and marginal distributions.\footnote{Our assumptions will guarantee the existence of $\rho$ (see Section~\ref{sec:assumptions}). Since $Y$ is a Polish space (because it is separable Hilbert) the $\rho(y|\mu_a)$ conditional 
    distribution ($y\in Y$, $\mu_a\in X$) is also well-defined \citep[Lemma~A.3.16, page~487]{steinwart08support}. \label{footnote:rho(|)}} 
    The regression function of $\rho$ with respect to the $(\mu_x,y)$ pair is denoted by
    \begin{align}
    f_{\rho}(\mu_a) &= \int_{Y}y\,\d \rho (y|\mu_a)\quad (\mu_a\in X) \label{eq:f_rho}
    \end{align}
    and for $f\in L^{2}_{\rho_X}$ let $\left\| f \right\|_{\rho} = \sqrt{\left<f,f\right>_{\rho}} := \left\| f \right\|_{L^2_{\rho_X}} = \left[\int_X \left\|f(\mu_a)\right\|_Y^2 \d \rho_X(\mu_a)\right]^{\frac{1}{2}}$.
    Let us assume that the operator-valued kernel $K:X\times X\rightarrow \L(Y)$  is a Mercer kernel (that is $\H=\H(K)\subseteq C(X,Y)=\{X\rightarrow Y\text{ continuous functions}\}$), is bounded 
    ($\exists B_K<\infty$ such that $\left\|K(x,x)\right\|_{\L(Y)}\le B_K$), and is a compact operator for all $x\in X$. These requirements will be guaranteed by our assumptions, see Section~\ref{proof:flambda-frho-with-h}.
    In this case,  the inclusion $S_K^*$: $\H\hookrightarrow L^2_{\rho_X}$ is bounded, and its adjoint $S_K:L^2_{\rho_X} \rightarrow \H$ is given by
    \begin{align}
	(S_Kg)(\mu_u)=\int_{X}K(\mu_u,\mu_t)g(\mu_t)\d \rho_X(\mu_t).\label{eq:T:inprecise}
    \end{align}
    We further define $\tilde{T}$  as
    \begin{align}
    \tilde{T} = S_K^* S_K: L^2_{\rho_X}\rightarrow L^2_{\rho_X};\label{eq:def:Ttilde}
    \end{align}
    in other words, the result of 
    operation \eqref{eq:T:inprecise} belongs to $\H$, 
    which is embedded in $L^2_{\rho_X}$.
    $\tilde{T}$ is a compact, positive, self-adjoint operator \citep[Proposition~3]{carmeli10vector}, thus by the spectral theorem $\tilde{T}^s$ exists, where $s\ge 0$.
\end{itemize}

\subsection{Distribution Regression}\label{sec:DR problem} We now formally define the distribution regression task. Let us assume that $\M^+_1(\X)$ is endowed with $\S_1=\Bo(\tau_w)$, the weak-topology generated $\sigma$-algebra; thus $(\M^+_1(\X),\S_1)$ is a measurable space.
In the \emph{distribution regression} problem, we are given samples $\hat{\b{z}} = \{(\{x_{i,n}\}_{n=1}^{N_i},y_i)\}_{i=1}^l$ with $x_{i,1},\ldots, x_{i,N_i} \stackrel{i.i.d.} {\sim} x_i$ 
where $\b{z}=\{(x_i,y_i)\}_{i=1}^l$ with $x_i\in \M^+_1\left(\X\right)$ and $y_i\in Y$  drawn i.i.d.~from a joint meta distribution 
$\mathcal{M}$ defined on the measurable space $(\M^+_1(\X)\times Y, \S_1  \otimes \Bo(Y))$, the product space enriched with the product $\sigma$-algebra. Unlike in classical supervised learning problems, the problem at hand involves two levels of randomness, wherein first $\b{z}$ is drawn from $\mathcal{M}$, and then $\hat{\b{z}}$ is generated by sampling points from $x_i$ for all $i=1,\ldots,l$. The goal is to learn the relation between the random distribution $x$ and response $y$ based on 
the observed $\hat{\b{z}}$. For notational simplicity, we will assume that $N=N_i$ ($\forall i$).

As in the classical regression problem ($\R^d\rightarrow \R$), distribution regression can be tackled via kernel ridge regression  (using a squared loss as the discrepancy criterion). The 
kernel (say $\Kmeta$) is defined on $\M_1^+(\X)$, and the regressor is then modelled  by an element in the RKHS $\Hmeta=\Hmeta(\Kmeta)$ of functions mapping from $\M_1^+(\X)$ to $Y$. In this paper, we choose 
$\Kmeta(x,x')=K(\mu_x,\mu_{x'})$ where $x,x'\in \M_1^+(\X)$ and so that 
the function (in $\Hmeta$) to describe the $(x,y)$ random relation is constructed as a composition $f\circ \mu_x$, i.e.
\begin{align}
  \M^+_1\left(\X\right) \xrightarrow{\mu} X (\subseteq H=H(k)) \xrightarrow{f\in \H=\H(K)} Y. \label{eq:composition}
\end{align}
In other words, the distribution $x\in \M^+_1\left(\X\right)$ is first mapped to $X\subseteq H$ by the mean embedding $\mu$, and the result is composed with $f$, 
an element of the RKHS $\H$. 

Let the expected risk for a $\tilde{f}:X\rightarrow Y$ (measurable) function be defined as 
      \begin{align*}
	  \Eo\big[\tilde{f}\big] &= \E_{(x,y)\sim\mathcal{M}} \big\|\tilde{f}(\mu_x) - y\big\|_Y^2,			          
      \end{align*}
which is minimized by the $f_{\rho}$ regression function. The classical regularization approach is to optimize
\begin{align}
	f_{\b{z}}^{\lambda} &= \argmin_{f\in \H}\frac{1}{l}\sum_{i=1}^l\left\|f(\mu_{x_i})-y_i\right\|_Y^2 + \lambda \left\|f\right\|_{\H}^2\label{eq:obj0}
\end{align}
 instead of $\Eo$, based on samples $\b{z}$. Since $\b{z}$ is not available, we consider the objective function defined by 
the observable quantity $\hat{\b{z}}$,
\begin{align}
    f_{\hat{\b{z}}}^{\lambda} &= \argmin_{f\in \H}\frac{1}{l}\sum_{i=1}^l\left\|f(\mu_{\hat{x}_i}) - y_i\right\|_Y^2 + \lambda \left\|f\right\|_{\H}^2, \label{eq:obj}
\end{align}
where $\hat{x}_i=\frac{1}{N}\sum_{n=1}^N \delta_{x_{i,n}}$ is the empirical distribution determined by $\left\{x_{i,n}\right\}_{i=1}^N$. The ridge regression objective function has an analytical solution: given training samples $\hat{\b{z}}$, the prediction for a new $t$ test distribution is
\begin{align}
    (f_{\hat{\b{z}}}^{\lambda} \circ \mu) (t) &= \b{k} (\b{K}+l \lambda I)^{-1} [y_1;\ldots;y_l], \label{eq:MERR:analytical-solution}
\end{align}
where $\b{k} =\left[K(\mu_{\hat{x}_1},\mu_t), \ldots, K(\mu_{\hat{x}_l},\mu_t)\right]\in\L(Y)^{1\times l}$, $\b{K} =  [K(\mu_{\hat{x}_i},\mu_{\hat{x}_j})]\in \L(Y)^{l\times l}$, $[y_1;\ldots;y_l]\in Y^l$.

\begin{remark}~
\begin{itemize}
  \item It is important to note that the algorithm has access to the sample points \emph{only via} their \emph{mean embeddings} $\{\mu_{\hat{x}_i}\}_{i=1}^l$ in Eq.~\eqref{eq:obj}.
  \item There is a \emph{two-stage sampling difficulty} to tackle: The transition from $f_{\rho}$ to $f_{\b{z}}^{\lambda}$ represents the fact that we have only $l$ distribution samples ($\b{z}$); 
	the transition from $f_{\b{z}}^{\lambda}$ to $f_{\hat{\b{z}}}^{\lambda}$ means that the $x_i$ distributions can be accessed only via samples 
	($\hat{\b{z}}$).
  \item While  ridge regression can be performed using the kernel $\Kmeta$, the two-stage sampling makes it difficult to work with arbitrary $\Kmeta$. By contrast, our 
	choice of $\Kmeta(x,x')=K(\mu_x,\mu_{x'})$ enables us to handle the two-stage sampling by estimating $\mu_x$ with an empirical estimator, and using it in the algorithm as shown above.
  \item In case of scalar output ($Y=\R$), $\L(Y)=\L(\R)=\R$ and \eqref{eq:MERR:analytical-solution} is a standard linear equation with 
      $\b{K}\in \R^{l\times l}$, $\b{k}\in\R^{1\times l}$. More generally, if $Y=\R^d$, then $\L(Y)=\L(\R^d)=\R^{d\times d}$ and \eqref{eq:MERR:analytical-solution} is still a finite-dimensional 
      linear equation with $\b{K}\in \R^{(dl)\times (dl)}$ and $\b{k}\in\R^{d\times (dl)}$.
  \item \label{remarks:more-abstract-formulation} One could also formulate the problem (and get guarantees) for more abstract $X\subseteq H \rightarrow Y$ regression tasks 
  [see Eq.~\eqref{eq:composition}] on a convex set $X$ with $H$ and $Y$ being general, separable Hilbert spaces.
	Since distribution regression is probably the most accessible example where two-stage sampling appears, and in order to keep the 
	presentation simple, we do not consider such extended formulations in this work.
\end{itemize}
\end{remark}

Our main goals in this paper are as follows: first, to analyse the excess risk 
	\begin{align*}
	      \EoS \big(f^{\lambda}_{\hat{\b{z}}},f_{\rho}\big) := \Eo[f^{\lambda}_{\hat{\b{z}}}] -  \Eo[f_{\rho}],	 
	\end{align*}
	both when $f_{\rho}\in \H$ (the well-specified case) and $f_{\rho}\in L^2_{\rho_X}\backslash\H$ (the misspecified case); second, to establish
	consistency ($\EoS \left(f^{\lambda}_{\hat{\b{z}}},f_{\rho}\right)\rightarrow 0$, or in the misspecified case $\EoS \left(f^{\lambda}_{\hat{\b{z}}},f_{\rho}\right) - D_{\H}^2\rightarrow 0$, where
        $D_{\H}^2: = \inf_{q\in \H} \left\|f_{\rho}-S_K^*q\right\|_{\rho}^2$ is the approximation error of $f_{\rho}$ by a function in $\H$); and third, to derive an 
exact computational-statistical efficiency trade-off as a function of the $(l,N,\lambda)$ triplet, and of the difficulty of the problem.

\section{Assumptions} \label{sec:assumptions}
In this section, we detail our assumptions on the $(\X,Y,k,K)$ quartet. 
Our analysis for the well-specified case uses existing ridge regression results \citep{caponnetto07optimal} 
focusing on problem \eqref{eq:obj0} where only a single-stage sampling is present, hence we have to verify the associated conditions. 
Though we make use of these results, the analysis still remains challenging; the available bounds can moderately shorten our proof.
We must take particular care in verifying that \citet{caponnetto07optimal}'s conditions are met, since they need to hold for the space of \emph{mean embeddings of the distributions} ($X=\mu\left(\M_1^+(\X)\right)$), whose properties as a function of $\X$ and $H$ must themselves be established.\vspace{1mm}

\noindent Our \tb{assumptions} are as follows:
  \begin{enumerate}[topsep=0cm,partopsep=0cm,parsep=0cm,itemsep=0cm]
	  \item $(\X,\tau)$ is a separable, topological space. 
	  \item $Y$ is a separable Hilbert space.
	  \item $k$ is bounded, in other words $\exists B_k < \infty$ such that $\sup_{u\in\X}k(u,u)\le B_k$, and continuous.
	  \item The $\{K_{\mu_a}\}_{\mu_a\in X}$ operator family is uniformly bounded in Hilbert-Schmidt norm and H{\"o}lder continuous in operator norm. Formally, $\exists B_K < \infty$ such that
		    \begin{align}
		      \left\|K_{\mu_a}\right\|_{\L_2(Y,\H)}^2 = Tr \left(K_{\mu_a}^* K_{\mu_a}\right) \le B_K, \quad (\forall \mu_a\in X), \label{eq:bounded kernel}
		    \end{align}
		    and $\exists L>0$, $h\in (0,1]$ such that the mapping $K_{(\cdot)}: X \rightarrow \L(Y,\H)$ is H\"{o}lder continuous:
			\begin{align}
			      \left\|K_{\mu_a} - K_{\mu_b}\right\|_{\L(Y,\H)} &\le L \left\|\mu_a - \mu_b\right\|_H^h,\quad \forall (\mu_a,\mu_b)\in X\times X.\label{eq:K:Lip}
			\end{align}
	 \item $y$ is bounded: $\exists C<\infty$ such that $\left\|y\right\|_Y\le C$ almost surely. 
  \end{enumerate}
These requirements hold under mild conditions:  in Section~\ref{sec:assumptions-discussion}, we provide insight into the consequences of our assumptions, with several concrete illustrations (e.g. regression with set- and RBF-type kernels).
 
\section{Error Bounds, Consistency \& Computational-Statistical Efficiency Trade-off} \label{sec:convergence analysis}
In this section, we present our analysis of the consistency of the mean embedding based ridge regression (MERR) method.

Given the estimator ($f^{\lambda}_{\hat{\b{z}}}$) in Eq.~\eqref{eq:obj},
we
derive finite-sample high probability upper bounds (see Theorems~\ref{theo1} and \ref{theo:L2}) for the excess risk $\EoS\left(f^{\lambda}_{\hat{\b{z}}},f_{\rho}\right)$, and 
in the misspecified setting, for the excess risk compared to the best attainable value from $\H$, i.e., $\EoS\left(f^{\lambda}_{\hat{\b{z}}},f_{\rho}\right)-D_{\H}^2$.
We  illustrate the bounds for particular classes of prior distributions, and work through special cases
 to obtain consistency conditions and computational-statistical efficiency trade-offs (see Theorems~\ref{conseq:excess-rate},  \ref{conseq:L2rate:Assumption2b} and the 3rd bullet of Remark~\ref{remark:misspec}).
The main challenge  is how to turn the convergence rates of the mean embeddings into those for an error $\EoS$ of the predictor. 
Although the main ideas of the proofs can be summarized relatively briefly, the full details are more demanding. High-level ideas with the sketches of the proofs and the 
obtained results are presented in Section~\ref{sec:results:A1} (well-specified case) and Section~\ref{sec:results:A2} (misspecified case). The derivations of some technical details of
Theorems~\ref{theo1} and~\ref{theo:L2} are available in Section~\ref{sec:proofs}.

\subsection{Results for the Well-specified Case}\label{sec:results:A1}
We first focus on the well-specified case ($f_{\rho}\in \H$) and present our first main result. We derive a high probability upper bound for the excess risk 
$\EoS\left(f^{\lambda}_{\hat{\b{z}}},f_{\rho}\right)$ of the MERR
method (Theorem~\ref{theo1}). The upper bound is instantiated for a general class of prior distributions (Theorem~\ref{conseq:excess-rate}), which leads to a 
simple computational-statistical efficiency description (Theorem~\ref{conseq:conv-rate}); this shows (among others) conditions when the MERR technique is able to achieve the \emph{one-stage} sampled minimax optimal rate. We first give a high-level sketch of our convergence analysis and an intuitive interpretation of the results. 
An outline of the main proof ideas is given below, with technical details in Section~\ref{sec:proofs}.

Let us define $\b{x} = \{x_i\}_{i=1}^l$ and $\hat{\b{x}} = \{\{x_{i,n}\}_{n=1}^{N}\}_{i=1}^l$ as the `x-part' of $\b{z}$ and $\hat{\b{z}}$, respectively.
One can  express $f_{\b{z}}^{\lambda}$ [Eq.~\eqref{eq:obj0}] \citep{caponnetto07optimal}, and similarly $f_{\hat{\b{z}}}^{\lambda}$ [Eq.~\eqref{eq:obj}], as
\begin{align}
    f_{\b{z}}^{\lambda} &= (T_{\b{x}}+\lambda)^{-1}g_{\b{z}},& T_{\b{x}}& = \frac{1}{l}\sum_{i=1}^lT_{\mu_{x_i}}, & g_{\b{z}} &= \frac{1}{l}\sum_{i=1}^l K_{\mu_{x_i}}y_i \label{eq:f_zlambda},\\
    f_{\hat{\b{z}}}^{\lambda} &= (T_{\hat{\b{x}}}+\lambda)^{-1}g_{\hat{\b{z}}},& T_{\hat{\b{x}}} &= \frac{1}{l}\sum_{i=1}^lT_{\mu_{\hat{x}_i}}, & g_{\hat{\b{z}}}& = \frac{1}{l}\sum_{i=1}^l K_{\mu_{\hat{x}_i}} y_i \label{eq:f_hatz^lambda},
\end{align}
where  $T_{\mu_a} = K_{\mu_a}K_{\mu_a}^*\in \L(\H)$ ($\mu_a\in X$), $T_{\b{x}}, T_{\hat{\b{x}}}:\H\rightarrow \H$, $g_{\b{z}}, g_{\hat{\b{z}}}\in\H$. By these explicit expressions, one can decompose
the excess risk into 5 terms \citep[Section~A.1.8]{szabo15twostage}:
\begin{align*}
	      \EoS\big(f^{\lambda}_{\hat{\b{z}}},f_{\rho}\big) &= \Eo\big[f^{\lambda}_{\hat{\b{z}}}\big] - \Eo\left[f_{\rho}\right] \le 5 \left[S_{-1} + S_0 + \A(\lambda) + S_1 + S_2 \right],
\end{align*}
where
\begin{align}
	      S_{-1} &= S_{-1}(\lambda,\b{z},\hat{\b{z}}) = \|\sqrt{T} (T_{\hat{\b{x}}}+\lambda I)^{-1}(g_{\hat{\b{z}}} - g_{\b{z}})\|_{\H}^2,\label{eq:Assumption1:S_-1:def}\\ 
	      S_0 &= S_0(\lambda,\b{z},\hat{\b{z}}) = \|\sqrt{T} (T_{\hat{\b{x}}}+\lambda I)^{-1}(T_{\b{x}}-T_{\hat{\b{x}}}) f^{\lambda}_{\b{z}}\|_{\H}^2,\label{eq:Assumption1:S_0:def}\\
	      \A(\lambda) &= \|\sqrt{T}(f^{\lambda}-f_{\rho})\|_{\H}^2, \hspace*{0.4cm}
	      S_1 = S_1(\lambda,\b{z}) = \|\sqrt{T}(T_{\b{x}}+\lambda I)^{-1}(g_{\b{z}}-T_{\b{x}}f_{\rho})\|_{\H}^2,\nonumber\\
	      S_2&= S_2(\lambda,\b{z}) = \|\sqrt{T}(T_{\b{x}}+\lambda I)^{-1}(T-T_{\b{x}})(f^{\lambda}-f_{\rho})\|_{\H}^2,\nonumber\\
  f^{\lambda} &= \argmin_{f\in\H}(\Eo[f]+\lambda\left\|f\right\|_{\H}^2), \hspace{0.2cm} T = \int_{X}T_{\mu_a}\d \rho_X(\mu_a)=S_K S_K^*:\H\rightarrow \H \label{eq:flambda:def}.
\end{align}
Three of the terms ($S_1$, $S_2$, $\A(\lambda)$) are identical to the terms in \citet{caponnetto07optimal}, hence the earlier bounds can be applied.
The two new terms ($S_{-1}$, $S_0$) resulting from  two-stage sampling will be upper bounded by making use of the convergence of the empirical mean embeddings.
These bounds will lead to the following results:
	      \begin{theorem}[Finite-sample excess risk bounds; well-specified case]\label{theo1}
		    Let $K_{(\cdot)}: X\rightarrow \L(Y,\H)$ be H{\"o}lder continuous with constants $L$, $h$.
		    Let $l\in\Z^+$, $N\in\Z^+$, $0<\lambda$,  $0<\eta<1$, $0<\delta$, $C_{\eta}=32\log^2(6/\eta)$, $\left\|y\right\|_Y\le C$ (a.s.) and $\A(\lambda)$ be the residual 
		    as defined above. Define $M =2(C + \left\|f_{\rho}\right\|_{\H}\sqrt{B_K})$, $\Sigma = \frac{M}{2}$, $T$ as in \eqref{eq:flambda:def}, 
		    $\B(\lambda)=\|f^{\lambda} - f_{\rho}\|_{\H}^2$ as the reconstruction error, and $\N(\lambda) = Tr [(T+\lambda I)^{-1}T]$ as the effective dimension.
		    Then with probability at least $1-\eta-e^{-\delta}$, the excess risk can be upper bounded as
		    \begin{eqnarray*}
			\lefteqn{\EoS\left(f^{\lambda}_{\hat{\b{z}}},f_{\rho}\right)
			      \le 5\Bigg\{ \frac{4L^2\left(1+\sqrt{\log(l)+\delta}\right)^{2h}(2B_k)^h}{\lambda N^h}\left[C^2+4B_K\times \phantom{\left(\log^2\left(\frac{6}{\eta}\right) \left\{\frac{64}{\lambda} \left[\frac{M^2 B_K}{l^2\lambda} + \frac{\Sigma^2\N(\lambda)}{l} \right]
				+ \frac{24}{\lambda^2} \left[ \frac{4B_K^2\B(\lambda)}{l^2} + \frac{B_K \A(\lambda)}{l}\right]\right\}  + \B(\lambda) + \left\| f_{\rho}\right\|_{\H}^2\right)}\right.}\\
			      && \left.\hspace*{-0.5cm} \times \left(\log^2\left(\frac{6}{\eta}\right) \left\{\frac{64}{\lambda} \left[\frac{M^2 B_K}{l^2\lambda} + \frac{\Sigma^2\N(\lambda)}{l} \right]
				+ \frac{24}{\lambda^2} \left[ \frac{4B_K^2\B(\lambda)}{l^2} + \frac{B_K \A(\lambda)}{l}\right]\right\}  + \B(\lambda) + \left\| f_{\rho}\right\|_{\H}^2\right)\right]\nonumber\\
			      &&\hspace*{4.7cm}  + \A(\lambda) + C_{\eta} \left[\frac{B_K^2\B(\lambda)}{l^2\lambda} + \frac{B_K \A(\lambda)}{4l\lambda} + \frac{B_K M^2}{l^2\lambda} + \frac{\Sigma^2 \N(\lambda)}{l} \right] \Bigg\}
		    \end{eqnarray*}
		    if
			    $l\ge 2C_{\eta}B_K\N(\lambda) / \lambda$,  $\lambda \le \left\|T\right\|_{\L(\H)}$ and 
			    $N \ge \big(1+\sqrt{\log(l)+\delta}\big)^{2} 2^{\frac{h+6}{h}}B_k (B_K)^{\frac{1}{h}} L^{\frac{2}{h}} / \lambda^{\frac{2}{h}}$.
	      \end{theorem}
	      Below we specialize our excess risk bound for a general prior class, which captures the difficulty of the regression problem as defined in \citet{caponnetto07optimal}. This 
    	      $\P(b,c)$ class is described by two parameters $b$ and $c$: larger $b$ means faster decay of the eigenvalues of the covariance operator $T$ 
	      [in  Eq. \eqref{eq:flambda:def}], hence smaller effective input dimension; larger $c$  corresponds  to a smoother regression function. Formally:\vspace{1mm}
	
	      \noindent\tb{Definition of the $\P(b,c)$ class:} Let us fix the positive constants $R$, $\alpha$, $\beta$. Then given $1<b$, $c\in(1,2]$, the $\P(b,c)$ class is
	      the set of probability distributions $\rho$ on $Z=X\times Y$  such that 
	      \begin{enumerate}[topsep=0cm,partopsep=0cm,parsep=0cm,itemsep=0cm]
		  \item a range space assumption is satisfied: $\exists g\in\H$ s.t.\ $f_{\rho}=T^{\frac{c-1}{2}}g$ with $\left\|g\right\|_{\H}^2\le R$,
		  \item in the spectral decomposition of $T=\sum_{n=1}^{\infty}\lambda_n\left<\cdot,e_n\right>_{\H}e_n$, where $(e_n)_{n=1}^{\infty}$ is a basis of $\Ker(T)^{\perp}$, the eigenvalues of 
			$T$ satisfy $\alpha \le n^b \lambda_n\le \beta \quad (\forall n\ge 1)$.%\vspace{1mm}
	      \end{enumerate}
	      \begin{remark} We make few remarks about the $\mathcal{P}(b,c)$ class:\vspace{-1mm}
   	      \begin{itemize}
		\item \emph{Range space assumption on  $f_{\rho}$:} \label{remark:P(b,c)-interpretation} The smoothness of $f_{\rho}$ is expressed as a range space assumption, which is slightly different from the standard smoothness conditions 
		      appearing in non-parametric function estimation. By the spectral decomposition of $T$ given above [$\lambda_1\ge \lambda_2 \ge \ldots >0, \lim_{n\rightarrow \infty}\lambda_n = 0$],
		       $T^r u = \sum_{n=1}^{\infty} (\lambda_n)^r \left<u,e_n\right>_{\H}e_n \quad (r=\frac{c-1}{2}\ge 0, u\in \H)$ and
		      \begin{align}
		       \Im(T^r) = \Big\{\sum\nolimits_{n=1}^{\infty} c_n e_n: \sum\nolimits_{n=1}^{\infty} c_n^2\lambda_n^{-2r}<\infty\Big\}. \label{eq:range}
		      \end{align}
		      Specifically, in the limit as $r\rightarrow 0$, we obtain $f_{\rho}\in \Im(T^0)=\Im(I) = \H$ (no constraint); larger values of 
		      $r$ give rise to faster decay of the $(c_n)_{n=1}^{\infty}$ Fourier coefficients. This is the concrete  meaning of $f_{\rho}\in \Im(T^r)$.
		\item \emph{Spectral decay condition:} We can provide a simple illustration of when the spectral decay conditions hold, 
		      in the event that the distributions are normal with means $m_i$ and identical variance ($x_i = N(m_i,\sigma^2I$)).
		      When Gaussian kernels ($k$) are used with linear $K$, 
		      then $K(\mu_{x_i},\mu_{x_j}) = e^{-c\left\|m_i-m_j\right\|^2}$ \citep[Table 1, line 2]{muandet12learning}  
		      (Gaussian, with arguments equal to the difference in means). Thus, this Gram matrix will correspond to the Gram matrix using a Gaussian kernel between points 
		      $m_i$. The spectral decay of the Gram matrix will correspond to that of the Gaussian kernel, with points drawn from the meta-distribution over the $m_i$. 
		      Thus, the source conditions are analysed in the same manner as for Gaussian Gram matrices:  see e.g. \citet{steinwart08support} for a discussion of these spectral decay properties.\vspace{1mm}
   	      \end{itemize}
   	      \end{remark}
		In the $\P(b,c)$ family, the behaviour of $\A(\lambda)$, $\B(\lambda)$ and $\N(\lambda)$ is known: $\A(\lambda)\le R\lambda^c$, $\B(\lambda)\le R\lambda^{c-1}$, $\N(\lambda)\le \beta \frac{b}{b-1}\lambda^{-\frac{1}{b}}$. Specializing Theorem~\ref{theo1} and retaining its assumptions, we get:
	      \begin{theorem}[Finite-sample excess risk bound for $\rho\in \P(b,c)$]\label{conseq:excess-rate} Suppose the conditions in Theorem~\ref{theo1} hold. Let $\rho\in \P(b,c)$, where $1<b$ and $c\in (1,2]$. Then
	      \begin{align*}
		  \EoS\left(f^{\lambda}_{\hat{\b{z}}},f_{\rho}\right)
		  &\le 5 \Bigg\{ \frac{4L^2 \left(1+\sqrt{\log(l)+\delta}\right)^{2h} (2B_k)^h}{\lambda N^h} \Bigg[  C^2 + 4B_K\times  \\
		  & \left.\hspace*{-1.3cm} \times \left(C_{\eta} \left\{\frac{2}{\lambda} \left[\frac{M^2 B_K}{l^2\lambda} + \frac{\Sigma^2\beta b}{(b-1)l\lambda^{\frac{1}{b}}} \right]
				+ \frac{3}{4\lambda^2} \left[ \frac{4B_K^2R\lambda^{c-1}}{l^2} + \frac{B_K R\lambda^c}{l}\right]\right\}  + R\lambda^{c-1} + \left\| f_{\rho}\right\|_{\H}^2\right)\right]\\
		  &\qquad  + R\lambda^c + C_{\eta} \left[\frac{B_K^2R\lambda^{c-2}}{l^2}+\frac{B_K R\lambda^{c-1}}{4l}+\frac{B_K M^2}{l^2\lambda}+\frac{\Sigma^2\beta b}{(b-1)l\lambda^{\frac{1}{b}}}\right]\Bigg\}.
	      \end{align*}	      
	      \end{theorem}
	      Discarding the constants in Theorem~\ref{conseq:excess-rate}, the study of convergence of the excess risk $\EoS(f^{\lambda}_{\hat{\b{z}}},f_{\rho})$ to $0$ boils down to finding $N$ and $\lambda$ (as a function of $l$) where $N\rightarrow\infty$, $\lambda\rightarrow 0$ and
	      \begin{align}
			    r(l,N,\lambda) &= \frac{\log^h(l)}{N^h\lambda}\left(\frac{1}{\lambda^2l^2} + 1+ \frac{1}{l\lambda^{1+\frac{1}{b}}}\right)+\lambda^c + \frac{1}{l^2\lambda} + \frac{1}{l\lambda^{\frac{1}{b}}}\rightarrow 0, 
		    \text{ s.t. } l \lambda^{\frac{b+1}{b}}\ge 1, \frac{\log(l)}{\lambda^{\frac{2}{h}}} \le N\label{eq:f1}
	      \end{align}
	      as $l\rightarrow\infty$. Let us choose $N=l^{\frac{a}{h}}\log(l)$; in this case Eq.~\eqref{eq:f1} reduces to
	      \begin{align}
	      r(l,\lambda) 
		&= \frac{1}{l^{2+a}\lambda^3} + \frac{1}{l^a\lambda} + \frac{1}{l^{a+1}\lambda^{2+\frac{1}{b}}} + \lambda^c + \frac{1}{l^2\lambda} + \frac{1}{l\lambda^{\frac{1}{b}}}\rightarrow 0, 
		    \text{ s.t. } l \lambda^{\frac{b+1}{b}}\ge 1,\hspace{0.1cm} l^a\lambda^2\ge 1 \label{eq:f1b}.
	      \end{align}
	      One can assume that $a>0$, otherwise $r(l,\lambda)\rightarrow 0$ fails to hold; in other words, $N$ should grow faster than $\log(l)$.  
	      Matching the `bias' ($\lambda^s$) and `variance' (other) terms in
	      $r(l,\lambda)$ to choose  $\lambda$, and guaranteeing that the matched terms dominate and the constraints in Eq.~\eqref{eq:f1b} hold, one gets the following simple 
		description for the computational-statistical efficiency 
		trade-off:\footnote{The derivations are available in the supplement.\label{footnote:rate-derivations}}

	    \begin{theorem}\emph{\tb{(Computational-statistical efficiency trade-off; well-specified case; $\rho\in \P(b,c)$)}} \label{conseq:conv-rate} Suppose the conditions in Theorem~\ref{theo1} hold. Let $\rho\in \P(b,c)$ and $N=l^{\frac{a}{h}}\log(l)$, where $0<a$, $1<b$, $c\in (1,2]$. If
	      \begin{itemize}
		\item $a \le \frac{b(c+1)}{bc+1}$, then $\EoS\left(f^{\lambda}_{\hat{\b{z}}},f_{\rho}\right)=\O_{p}\left(l^{-\frac{ac}{c+1}}\right)$ with $\lambda = l^{-\frac{a}{c+1}}$,
		\item $a\ge \frac{b(c+1)}{bc+1}$ then $\EoS\left(f^{\lambda}_{\hat{\b{z}}},f_{\rho}\right)=\O_{p}\left(l^{-\frac{bc}{bc+1}}\right)$ with $\lambda = l^{-\frac{b}{bc+1}}$.
	      \end{itemize}
	      \end{theorem}
	      \begin{remark}
	      Theorem~\ref{conseq:conv-rate} formulates an exact computational-statistical efficiency trade-off for the choice of the bag size ($N$) as a  function of the number of distributions ($l$) and problem difficulty ($b$, $c$).\vspace{1mm}
	      \begin{itemize}
		  \item \emph{$a$-dependence:} A smaller bag size (smaller $a$; $N=l^{\frac{a}{h}}\log(l)$) means computational savings, but reduced statistical efficiency. It is not worth increasing $a$ above $\frac{b(c+1)}{bc+1}$ since from that point the rate 
			becomes $r(l)=l^{-\frac{bc}{bc+1}}$; remarkably, this rate is minimax in the one-stage sampled setup \citep{caponnetto07optimal}. The sensible choice $a= \frac{b(c+1)}{bc+1}<2$ means that the 
			one-stage sampled minimax rate can be achieved in the two-stage sampled setting with 
			bag size $N$ sub-quadratic in $l$.
		  \item \emph{$h$-dependence:} In accord with our `smoothness' assumptions it is rewarding to use smoother $K$ kernels (larger $h\in (0,1]$) since this reduces the bag size [$N=l^{\frac{a}{h}}\log(l)$].
		  \item \emph{$c$-dependence:} The strictly decreasing property of $c\mapsto  \frac{b(c+1)}{bc+1}$ implies that for `smoother' problems (larger $c$) fewer samples ($N$) are sufficient.
	      \end{itemize}
	      \end{remark}

Below we elaborate on the sketched high-level idea and prove Theorem~\ref{theo1}.\vspace{1mm}\\
\noindent\tb{Proof of Theorem~\ref{theo1}} (detailed derivations of each step can be found in Section~\ref{sec:proofs:well-specified})\vspace{1mm}
\begin{enumerate}[labelindent=0cm,leftmargin=*]
  \item \tb{Decomposition of the excess risk}:
 	We have the following upper bound for the excess risk
	\begin{align}
            \EoS\big(f^{\lambda}_{\hat{\b{z}}},f_{\rho}\big) &= \Eo\big[f^{\lambda}_{\hat{\b{z}}}\big] - \Eo\left[f_{\rho}\right] \le 5 \left[S_{-1} + S_0 + \A(\lambda) + S_1 + S_2 \right]. \label{eq:E:5terms}
	\end{align}
  \item \tb{It is sufficient to upper bound $S_{-1}$ and $S_0$}: \citet{caponnetto07optimal} have shown that for $\forall \eta>0$ if
   	      $l\ge 2C_{\eta}B_K\N(\lambda) / \lambda$, $\lambda \le \|T\|_{\L(\H)}$, 
	       then $\Pr(\bm{\Theta}(\lambda,\b{z})\le 1/2)\ge 1-\eta/3$, where
	      \begin{align}
		  \bm{\Theta}(\lambda,\b{z}) = \left\|(T - T_{\b{x}}) (T+\lambda I)^{-1}\right\|_{\L(\H)}, \label{eq:theta}
	      \end{align}
	      using which upper bounds on $S_1$ and $S_2$ that hold with probability $1-\eta$ are obtained. It is known that $\A(\lambda)\le R\lambda^c$.
  \item \tb{Probabilistic bounds on $\| g_{\hat{\b{z}}} - g_{\b{z}}\|_{\H}^2$, $\| T_{\b{x}}-T_{\hat{\b{x}}}\|_{\L(\H)}^2$, $\|\sqrt{T} (T_{\hat{\b{x}}}+\lambda I)^{-1}\|_{\L(\H)}^2$, $\|f_{\b{z}}^{\lambda}\|_{\H}^2$}:
	      One can bound $S_{-1}$ and $S_0$ as $$S_{-1} \le \big\| \sqrt{T} (T_{\hat{\b{x}}}+\lambda I)^{-1}\big\|_{\L(\H)}^2 \|g_{\hat{\b{z}}} - g_{\b{z}} \|_{\H}^2$$ and $$S_0 \le  \big\| \sqrt{T} (T_{\hat{\b{x}}}+\lambda I)^{-1}\big\|_{\L(\H)}^2 \| T_{\b{x}}-T_{\hat{\b{x}}}\|_{\L(\H)}^2 \big\|f^{\lambda}_{\b{z}}\big\|^2_{\H}.$$ 
	      For the terms on the r.h.s., we derive upper bounds [for the definition of $\alpha$, see Eq.~\eqref{eq:emp-mean-emb-conv-rate:Assumption1}]
	      $$\left\|  g_{\hat{\b{z}}} - g_{\b{z}}  \right\|_{\H}^2  \le  L^2 C^2  \frac{\left(1+\sqrt{\alpha}\right)^{2h} (2B_k)^h}{N^h}, \quad \left\| \sqrt{T} (T_{\hat{\b{x}}}+\lambda I)^{-1}\right\|_{\L(\H)}  \le \frac{2}{\sqrt{\lambda}},$$
	      $$\left\|T_{\b{x}} - T_{\hat{\b{x}}}\right\|_{\L(\H)}^2 \le  \frac{\left(1+\sqrt{\alpha}\right)^{2h} 2^{h+2}(B_k)^{h}B_K L^2}{N^h},$$
and
		  \begin{align}
		  \left\| f_{\b{z}}^{\lambda} \right\|_{\H}^2 &\le 6\left(\frac{16}{\lambda} \log^2\left(\frac{6}{\eta}\right)\left[\frac{M^2 B_K}{l^2\lambda} + \frac{\Sigma^2\N(\lambda)}{l} \right]\right.\label{f:zlambda:improved-bound}\\
		  &\quad \hspace*{0.6cm}\left. + \frac{4}{\lambda^2}\log^2\left(\frac{6}{\eta}\right) \left[ \frac{4B_K^2\B(\lambda)}{l^2} + \frac{B_K \A(\lambda)}{l}\right]  + \B(\lambda) + \left\| f_{\rho}\right\|_{\H}^2\right).\nonumber 
	      \end{align}
	      The bounds hold under the following conditions:
	      \begin{itemize}
		  \item $\| g_{\hat{\b{z}}} - g_{\b{z}}\|_{\H}^2$ (see Section~\ref{proof:gz}): if the empirical mean embeddings are close to their population counterparts, i.e.,
			    \begin{align}
				\left\|\mu_{x_i} - \mu_{\hat{x}_i}\right\|_H &\le \frac{(1+\sqrt{\alpha})\sqrt{2B_k}}{\sqrt{N}}, \quad (\forall i=1,\ldots,l).\label{eq:emp-mean-emb-conv-rate:Assumption1}
			  \end{align}
			  This event has probability $1-le^{-\alpha}$ over all $i=1,\ldots,l$ samples; see  \citep{altun06unifying} and \cite[Section~A.1.10]{szabo15twostage}.
		  \item $\| T_{\b{x}}-T_{\hat{\b{x}}}\|_{\L(\H)}^2$ (see Section~\ref{proof:Tx}): \eqref{eq:emp-mean-emb-conv-rate:Assumption1} is assumed.
		  \item\label{eq:bounds} $\|\sqrt{T} (T_{\hat{\b{x}}}+\lambda I)^{-1}\|_{\L(\H)}^2$ \citep[Section~A.1.11]{szabo15twostage}: \eqref{eq:emp-mean-emb-conv-rate:Assumption1},  $\bm{\Theta}(\lambda,\b{z})\le \frac{1}{2}$, and
			    \begin{align}
			      \frac{\left(1+\sqrt{\alpha}\right)^{2} 2^{\frac{h+6}{h}}B_k (B_K)^{\frac{1}{h}} L^{\frac{2}{h}}}{\lambda^{\frac{2}{h}}}  &\le N.			    	\label{eq:N}
			  \end{align}
		  \item $\|f_{\b{z}}^{\lambda}\|_{\H}^2$:  The bound is guaranteed to hold under the conditions of the bounds of $S_1$ and $S_2$.\textsuperscript{\ref{footnote:rate-derivations}}
	      \end{itemize}
  \item \tb{Union bound}: By applying an $\alpha=\log(l)+\delta$ reparameterization, and combining the received upper bounds with \citet{caponnetto07optimal}'s results 
      for $S_1$ and $S_2$, Theorem~\ref{theo1} follows (Section~\ref{proof:theo1:union bound}) with a union bound.\vspace{1mm}
\end{enumerate}

Finally, we note that existing results/ideas were used at two points to simplify our analysis: bounding $S_1$, $S_2$, $\Theta(\lambda,\b{z})$, $\left\|f_{\b{z}}^{\lambda}\right\|_{\H}^2$ \citep{caponnetto07optimal} and $\left\|\mu_{x_i}-\mu_{\hat{x}_i}\right\|_H$ \citep{altun06unifying}.\footnote{
We also corrected some constants in the previous works \citep{altun06unifying,caponnetto07optimal}.} 

\subsection{Results for the Misspecified Case}\label{sec:results:A2}
In this section, we focus on the misspecified case ($f_{\rho}\in L^2_{\rho_X}\backslash \H$) and present our second main result, which was inspired by the proof technique of \citet[Theorem 12]{sriperumbudur14density}. We derive a high probability upper bound for $\EoS\left(f^{\lambda}_{\hat{\b{z}}},f_{\rho}\right)$, i.e., the excess risk of the MERR
method (Theorem~\ref{theo:L2}) which gives rise to consistency results (3rd bullet of Remark~\ref{remark:misspec}) and precise computational-statistical efficiency trade-off (Theorem~\ref{conseq:L2rate:Assumption2b}). Theorem~\ref{theo:L2} consists of two finite-sample bounds:
\begin{enumerate}
 \item   The first, more general bound [Eq.~\eqref{eq:L2bound:explicit}] will be used to show consistency in the misspecified case (see the 3rd bullet of Remark~\ref{remark:misspec}), in other words that $\EoS \left(f^{\lambda}_{\hat{\b{z}}},f_{\rho}\right)$ can be driven to its smallest possible value determined by the ``richness'' of $\H$:
    \begin{align}
      D_{\H}^2: = \inf_{q\in \H} \left\|f_{\rho}-S_K^*q\right\|_{\rho}^2. \label{eq:DH}
    \end{align}
    The value of $D_{\H}$ equals the approximation error of $f_{\rho}$ by a function from $\H$. Specifically, if $\H$ [precisely $S_K^{*}(\H)=\{S_K^*q: q\in \H\}\subseteq L^2_{\rho_X}$] is \emph{dense} in $L^2_{\rho_X}$, 
	then $D_{\H}=0$.
  \item  The second, specialized result [Eq.~\eqref{eq:L2bound:explicit-b}] under additional smoothness assumptions on $f_{\rho}$ will give rise to  a precise 
	computational-statistical efficiency trade-off  in terms of the problem difficulty ($s$) and sample numbers ($l$, $N$); this result can be seen as the misspecified analogue of Theorem~\ref{conseq:conv-rate}.
\end{enumerate}

\noindent After stating our results, the main ideas of the proof follow; further technical details are available in Section~\ref{sec:proofs:misspecified}. Our main theorem for bounding the excess risk is as follows:
\begin{theorem}[Finite-sample excess risk bounds; misspecified case]\label{theo:L2}
    Let $l\in \Z^{+}$, $N\in\Z^{+}$, $0<\lambda$, $0<\eta<1$, $0<\delta$ and $C_{\eta}=\log\left(\frac{6}{\eta}\right)$. Assume that 
    $\left(\frac{12 B_K}{\lambda} C_{\eta}\right)^2 \le l$ and  $\big(1+\sqrt{\log(l)+\delta}\big)^{2} 2^{\frac{h+6}{h}}B_k (B_K)^{\frac{1}{h}} L^{\frac{2}{h}}  / \lambda^{\frac{2}{h}}\le N$.\vspace{1mm}
    \begin{enumerate}
    \item Then for arbitrary $q\in \H$ with probability at least $1-\eta-e^{-\delta}$
	    \begin{align}
		  \sqrt{\EoS \left(f^{\lambda}_{\hat{\b{z}}},f_{\rho}\right)} &\le
		  \frac{2LC\left(1+\sqrt{\log(l)+\delta}\right)^h(2B_k)^{\frac{h}{2}}}{\sqrt{\lambda}N^{\frac{h}{2}}} \left(1 + \frac{2\sqrt{B_K}}{\sqrt{\lambda}}\right) + \label{eq:L2bound:explicit}\\
		  &\hspace*{-1.3cm}\frac{2 C_\eta}{\sqrt{\lambda}} \left\{ \left( \frac{2C\sqrt{B}_K}{l}+ \frac{C\sqrt{B_K}}{\sqrt{l}} \right) + \left(\frac{2B_K}{l}+\frac{\sigma}{\sqrt{l}}\right) \frac{1}{\lambda}  \sqrt{\lambda \left\|f_{\rho}\right\|_{\rho} D_a(\lambda,q)} \right\} + D_a(\lambda,q)\nonumber,
	    \end{align}
	    where $D_a(\lambda,q) = \|f_{\rho}-S_K^*q\|_{\rho} + \max(1,\|T\|_{\L(\H)})\lambda^{\frac{1}{2}} \|q\|_{\H}$.\vspace{1mm}
    \item In addition, suppose $f_{\rho}\in \Im(\tilde{T}^s)$ for some $s>0$, where $\tilde{T}$ is defined in Eq.~\eqref{eq:def:Ttilde}. Then with probability at least $1-\eta-e^{-\delta}$, we have
	    \begin{align}
		\sqrt{\EoS \left(f^{\lambda}_{\hat{\b{z}}},f_{\rho}\right)} &\le \frac{2LC\left(1+\sqrt{\log(l)+\delta}\right)^h(2B_k)^{\frac{h}{2}}}{\sqrt{\lambda}N^{\frac{h}{2}}} \left(1 + \frac{2\sqrt{B_K}}{\sqrt{\lambda}}\right) +\nonumber\\
		  &\quad\hspace*{0cm} \frac{2 C_{\eta}}{\sqrt{\lambda}} \left\{ \left( \frac{2C\sqrt{B}_K}{l}+ \frac{C\sqrt{B_K}}{\sqrt{l}} \right) +  \left(\frac{2B_K}{l}+\frac{\sigma}{\sqrt{l}}\right) \frac{1}{\lambda} \times\right.\nonumber\\
		  &\quad\hspace{0.8cm} \left. \sqrt{\max\left(1,\left\|\tilde{T}\right\|^s_{\L\left(L^2_{\rho_X}\right)}\right)\lambda \left\|\tilde{T}^{-s}f_{\rho}\right\|_{\rho} D_b(\lambda,s)} \right\} + D_b(\lambda,s),\label{eq:L2bound:explicit-b}
	    \end{align}
	    where $D_b(\lambda,s) = \max(1,\|\tilde{T}\|_{\L(L^2_{\rho_X})}^{s-1}) \lambda^{\min(1,s)}\|\tilde{T}^{-s}f_{\rho}\|_{\rho}$.
    \end{enumerate}
\end{theorem}

\begin{remark} \label{remark:misspec} We give a short insight into the assumptions of Theorem~\ref{theo:L2}, followed by consequences of the theorem.
\begin{itemize}[labelindent=0cm,leftmargin=*,topsep=0cm,partopsep=0cm,parsep=0cm,itemsep=0cm]
 	  \item  \emph{Range space assumption on $f_{\rho}$:} \label{remark:rangespace-Ttilde} The range space assumption for the compact, positive, self-adjoint operator, $\tilde{T}=\tilde{T}(K): L^2_{\rho_X} \rightarrow L^2_{\rho_X}$ in the 2nd part of Theorem~\ref{theo:L2} can be interpreted similarly to 
 		that on $T$; see Eq.~\eqref{eq:range}. One can also prove alternative descriptions for $\Im(\tilde{T}^s)$ in terms of interpolation spaces \citep[Theorem~4.6, page 387]{steinwart12mercer}, or the decay of the 2-approximation error function, $A_2(\lambda)=\inf_{f\in \H(K)}\left(\lambda \left\|f\right\|_{\H(K)}^2+\Eo[f]-\Eo[f_{\rho}]\right)$ \citep{smale03estimating,steinwart09optimal}.
 	  \item $\sqrt{\EoS \left(f^{\lambda}_{\hat{\b{z}}},f_{\rho}\right)}$: \label{remark:square-root-of-excess-risk} Notice that in the bounds [\eqref{eq:L2bound:explicit}, \eqref{eq:L2bound:explicit-b}], instead of the excess risk, its square root appears; 
 	      this has technical reasons,  as it is easier to have the $D_a(\lambda,q)$ quantity (without multiplicative constants) appear on the r.h.s.\ of Eq.~\eqref{eq:L2bound:explicit} with this form.
	    \item \emph{Consistency in the misspecified case:} The consequence of Theorem~\ref{theo:L2}(1) is as follows.
	    Discarding the constants in Eq.~\eqref{eq:L2bound:explicit}, we obtain the upper bound (notice that the constant multiplier of $\left\|f_{\rho}-S_K^*q\right\|_{\rho}$ in the last term was \emph{one}):
	\begin{align*}
	    \sqrt{r(l,N,\lambda,q)}
		    &=\frac{\log^{\frac{h}{2}}(l)}{N^{\frac{h}{2}}\lambda} +  \frac{1}{\sqrt{l\lambda}} + \frac{\sqrt{\left\|f_{\rho}-S_K^*q\right\|_{\rho} + \sqrt{\lambda}\left\|q\right\|_{\H}}}{\lambda\sqrt{ l}} + \left\|f_{\rho}-S_K^*q\right\|_{\rho}
		    + \sqrt{\lambda} \left\|q\right\|_{\H}.
	\end{align*}
	By choosing $N=l^{1/h}\log l$, 
	  $\sqrt{r(l,\lambda)}$ is bounded by 
	  \begin{align*}
	      \inf_{q\in\mathcal{H}}\left\{\left\|f_{\rho}-S_K^*q\right\|_{\rho}+\frac{\sqrt{\left\|f_{\rho}-S_K^*q\right\|_{\rho}}}{\lambda\sqrt{ l}}+ \frac{\sqrt{\left\|q\right\|_{\H}}}{\lambda^{\frac{3}{4}}\sqrt{l}} +\sqrt{\lambda}\Vert q\Vert_\mathcal{H}\right\}+\O_{p}\left(\frac{1}{\sqrt{\lambda l}}\right).
	  \end{align*}
	  Our goal is to investigate the behavior of the bound as $l\rightarrow\infty$, $\lambda\rightarrow 0$ and $\lambda\sqrt{ l}\rightarrow \infty$. Define $K(\alpha,\beta,\gamma):=\inf_{q\in\mathcal{H}}\left\{\left\|f_{\rho}-S_K^*q\right\|_{\rho}+\alpha\sqrt{\left\|f_{\rho}-S_K^*q\right\|_{\rho}}+\beta\sqrt{\left\|q\right\|_{\H}}+\gamma\Vert q\Vert_\mathcal{H}\right\}$.
	  $K(\alpha,\beta,\gamma)$ is the pointwise infimum of affine functions, therefore it is upper semi-continuous and concave on $\R^3$ \citep[Lemmas~2.41 and 5.40 ]{aliprantis06infinite}; it is continuous on $\times_{i=1}^3\R_{>0}$ \citep[Theorem~2.35]{rockafellar08variational}. Moreover, by applying \citep[Corollary~2.37]{rockafellar08variational} it extends continuously to $\times_{i=1}^3\R_{\ge 0}$; specifically it is continuous at $(\alpha,\beta,\gamma)=\b{0}$. 
	   In other words, as $l\rightarrow\infty$, $\lambda\rightarrow 0$ and $\lambda\sqrt{l}\rightarrow \infty$,
	  $K\left(\frac{1}{\lambda\sqrt{ l}},\frac{1}{\lambda^{\frac{3}{4}}\sqrt{l}},\sqrt{\lambda}\right)\rightarrow D_\H$ and we get consistency in the misspecified case,\footnote{We have discarded the 
	  $\log(l)/ \lambda^{\frac{2}{h}}\le N$ constraint implied by the convergence of the first term in $\sqrt{r}$. \label{footnote:Ncondition:discard}} 
	  \begin{align*}
	      \sqrt{r(N,l,\lambda)} \rightarrow D_{\H}.
	  \end{align*}
\end{itemize}
\end{remark}

\begin{figure}
      \centering
      \includegraphics[width=5.5cm]{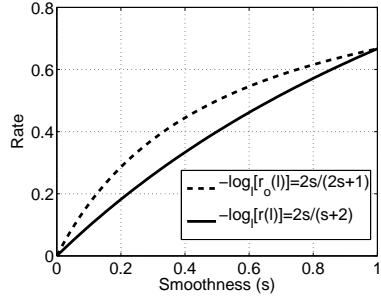}
      \caption{Comparison of the $r_o(l)=l^{-\frac{2s}{2s+1}}$ and $r(l)=l^{-\frac{2s}{s+2}}$ rates as function of the problem difficulty/smoothness ($s$).\label{fig:rates}}
  \end{figure}

	Discarding the constants in Eq.~\eqref{eq:L2bound:explicit-b} we get\textsuperscript{\ref{footnote:Ncondition:discard}}
	\begin{align}
	    \sqrt{r(l,N,\lambda)} 
			   &= \frac{\log^{\frac{h}{2}}(l)}{N^{\frac{h}{2}}\lambda} +  \frac{1}{\sqrt{l\lambda}} + \frac{\sqrt{\lambda^{\min(1,s)}}}{\lambda\sqrt{l}} + \lambda^{\min(1,s)} ,\text{ subject to } \frac{1}{\lambda^2}\le l. \label{eq:r:L2:rangespace}
	\end{align}
	Our goal is to drive $r(l,N,\lambda)$ to zero with a suitable choice of the $(l,N,\lambda)$ triplet under the stronger range space assumption.
	Since in Eq.~\eqref{eq:r:L2:rangespace} $\min(1,s)$ appears, one can assume without loss of generality that $s\in (0,1]$; consequently  
	$1-\frac{s}{2}\in \left[\frac{1}{2},1\right)$ and $\frac{1}{l^{\frac{1}{2}}\lambda^{\frac{1}{2}}} \le \frac{1}{\lambda^{1-\frac{s}{2}}l^{\frac{1}{2}}}$.
	Let us choose $N=l^{2a /h}\log(l)$; in this case using the previous dominance note, Eq.~\eqref{eq:r:L2:rangespace} reduces to the study of 
	\begin{align}
	    \sqrt{r(l,\lambda)} &= \frac{1}{l^a\lambda} +  \frac{1}{\lambda^{1-\frac{s}{2}}l^{\frac{1}{2}}} + \lambda^{s}\rightarrow 0,\text{ s.t. } l \lambda^2 \ge 1. \label{eq:r:L2:rangespace-temp2a}
	\end{align}
	One can assume that $a>0$, otherwise $r(l,\lambda)\rightarrow 0$ fails to hold: in other words, $N$ should grow faster than $\log(l)$. Matching the `bias' ($\lambda^s$) and `variance' (other) terms in
	$r(l,\lambda)$ to choose  $\lambda$, guaranteeing that the matched terms dominate and the constraint in Eq.~\eqref{eq:r:L2:rangespace-temp2a} hold, one can arrive at the following computational-statistical efficiency trade-off:\textsuperscript{\ref{footnote:rate-derivations}}
	\begin{theorem}\emph{\tb{(Computational-statistical efficiency trade-off; misspecified case, $f_{\rho}\in \Im(\tilde{T}^s)$)}}\label{conseq:L2rate:Assumption2b} 
	Suppose that $f_{\rho}\in \Im(\tilde{T}^s)$  and
	$N=l^{\frac{2a}{h}}\log(l)$,  where $s\in (0,1]$, $a>0$. If
	\begin{itemize}[labelindent=0.45cm,leftmargin=*,topsep=0.2cm,partopsep=0cm,parsep=0cm,itemsep=0cm]
	  \item $a\le \frac{s+1}{s+2}$, then $\EoS \left(f^{\lambda}_{\hat{\b{z}}},f_{\rho}\right)=\O_p\left(l^{-\frac{2sa}{s+1}}\right)$ with $\lambda = l^{-\frac{a}{s+1}}$, 
	  \item $a\ge \frac{s+1}{s+2}$, then $\EoS \left(f^{\lambda}_{\hat{\b{z}}},f_{\rho}\right)=\O_p\left(l^{-\frac{2s}{s+2}}\right)$ with $\lambda=l^{-\frac{1}{s+2}}$. 
	\end{itemize}
	\end{theorem}
	\begin{remark}
	Theorem~\ref{conseq:L2rate:Assumption2b} provides a complete computational-statistical efficiency trade-off description for the choice of the bag size $(N)$ as a number of the distributions ($l$).\vspace{1mm}
	  \begin{itemize}[labelindent=0.45cm,leftmargin=*,topsep=0cm,partopsep=0cm,parsep=0cm,itemsep=0cm]
	      \item \emph{$a$-dependence:} A smaller value of `$a$' (smaller bags $N=l^{2a /h}\log(l)$) leads to a computational advantage, but one looses in 
		    statistical efficiency. As `$a$' reaches $\frac{s+1}{s+2}$, the rate becomes $r(l)=l^{-\frac{2s}{s+2}}$ and one does not gain from further increasing the value of $a$. 
	            The sensible choice of $a = \frac{s+1}{s+2}\le \frac{2}{3}$ means that $N$ can again be \emph{sub-quadratic} ($2a<\frac{4}{3}<2$) in $l$. \vspace{1mm}
	      \item \emph{$h$-dependence:} By using smoother $K$ kernels (larger $h\in (0,1]$) one can reduce the size of the bags: $h\mapsto 2a/h$ is decreasing in $h$. This is compatible with our smoothness requirement on $f_{\rho}$.\vspace{1mm}
	      \item \emph{$s$-dependence:} ``Easier'' tasks (larger $s$) give rise to
			faster convergence. Indeed, in the $r(l)=l^{-\frac{2s}{s+2}}$ rate the $s\mapsto \frac{2s}{s+2}$ exponent is strictly increasing function of the problem difficulty ($s$). For example, 
			      for extremely non-smooth regression problems ($s\approx 0$) the convergence can be arbitrary slow ($\lim_{s\rightarrow 0} \frac{2s}{s+2} = 0$). In the smooth case ($s=1$)
			      $\lim_{s\rightarrow 1} \frac{2s}{s+2} = \frac{2}{3}$ and one can achieve the $r(l)=l^{-\frac{2}{3}}$ rate.\vspace{1mm}
	      \item We may compare our $r(l)=l^{-\frac{2s}{s+2}}$ result with the $r_o(l)=l^{-\frac{2s}{2s+1}}$ (one-stage sampled) rate \citep[$\beta/2:=s$, $q:=2$, $p:=1$ in Corollary~6]{steinwart09optimal}, which was shown to be asymptotically optimal on $Y=\R$ for continuous $k$ on compact metric $\X$. Steinwart {\em et al.}'s result is more general in terms of $q$ ($\left\|f\right\|_{\H}^q$ based regularization) and $p$ ($\left\|f\right\|_{\infty}\le C \left\|f\right\|_{\H}^p\left\|f\right\|_{\rho}^{1-p}$, $\forall f\in \H$; in our case $p=1$), although it imposes an \emph{additional} eigenvalue constraint [\cite[Eq.~(6)]{steinwart09optimal}] as well as $f_{\rho}\in \Im(\tilde{T}^s)$. Moreover, one can observe that $r_o(l)\le r(l)$ with a small gap, and that for $s\rightarrow 0$ and $s=1$, $r_o(l)=r(l)$; see Fig.~\ref{fig:rates}. We further remind the reader that our MERR analysis also holds for  \emph{separable Hilbert} output spaces $Y$,  \emph{separable topological} domains $\X$ enriched with a bounded, continuous kernel $k$, and that we handle the \emph{two-stage sampled} setting.
	\end{itemize}
\end{remark}

The main steps of the proof of Theorem~\ref{theo:L2} are as follows:\vspace{1mm}\\
\noindent\tb{Proof of Theorem~\ref{theo:L2}} (the details of the derivation are available in Section~\ref{sec:proofs:misspecified}) Steps~1-7 will be identical in both proofs,\footnote{Importantly, with a slight modification of the 
    more general, first part of Theorem~\ref{theo:L2}, one can get the specialized second setting of the theorem (see Step 8).\label{footnote:Th2:part1-2:diff}} and we present them jointly.
  \begin{enumerate}[labelindent=0cm,leftmargin=*]
      \item \tb{Decomposition of the excess risk}: By the triangle inequality, we have
	    \begin{align}
	      \sqrt{\EoS \left(f^{\lambda}_{\hat{\b{z}}},f_{\rho}\right)} &= \big\|S_K^* f^{\lambda}_{\hat{\b{z}}}- f_{\rho}\big\|_{\rho}   
								   \le \big\|S_K^* \big(f^{\lambda}_{\hat{\b{z}}} - f^{\lambda}_{\b{z}}\big)\big\|_{\rho} + \big\|S_K^* f^{\lambda}_{\b{z}} - f_{\rho}\big\|_{\rho}. \label{eq:L2}
	    \end{align}
      \item \tb{Bound on $\left\|S_K^* \left(f^{\lambda}_{\hat{\b{z}}} - f^{\lambda}_{\b{z}}\right)\right\|_{\rho}$}: 
	    Using\footnote{See for example \citet{devito06discretization} on page 88 with the $(\H,\G,A,T):= (\H,L^2_{\rho_X},S_K^*,T)$ choice.} the fact that 
	      \begin{align}
		    \left\|S_K^*h\right\|_\rho^2 &= \big\|\sqrt{T}h\big\|_{\H}^2 \quad (\forall h\in \H), \label{eq:S_K^*:norm preserving}
		  \end{align}
	     and the definitions of $S_{-1}$ and $S_0$ [see Eqs.~\eqref{eq:Assumption1:S_-1:def}-\eqref{eq:Assumption1:S_0:def}], we obtain
	    \begin{align}
	      \big\|S_K^* \big(f^{\lambda}_{\hat{\b{z}}} - f^{\lambda}_{\b{z}}\big)\big\|_{\rho} &= \big\|\sqrt{T} \big(f^{\lambda}_{\hat{\b{z}}} - f^{\lambda}_{\b{z}}\big)\big\|_{\H} \le \sqrt{S_{-1}} + \sqrt{S_0},\label{eq:L2:first-term}
	    \end{align}
	    through an application of triangle inequality. One can derive without a $\P(b,c)$ prior assumption  (Section~\ref{proof:S_-1+S_0}) the upper bound\footnote{See the remark at the end of Section~\ref{proof:S_-1+S_0}.}
	    \begin{align*}
	      \sqrt{S_{-1}} + \sqrt{S_0} &\le \frac{2LC(1+\sqrt{\alpha})^h(2B_k)^{\frac{h}{2}}}{\sqrt{\lambda}N^{\frac{h}{2}}} \left[1 + \frac{2\sqrt{B_K}}{\sqrt{\lambda}}\right]
	    \end{align*}
	     for the r.h.s.\ of Eq.~\eqref{eq:L2:first-term} under the conditions
	    that $\bm{\Theta}(\lambda,\b{z}) \le \frac{1}{2}$ (which holds with probability $1-\eta$ if
		$\left[12 B_K \log(2/\eta)/\lambda\right]^2 \le l$), and that Eqs.~\eqref{eq:emp-mean-emb-conv-rate:Assumption1}-\eqref{eq:N} hold.

      \item \tb{Decomposition of $\big\|S_K^* f^{\lambda}_{\b{z}} - f_{\rho}\big\|_{\rho}$}: By the triangle inequality and Eq.~\eqref{eq:S_K^*:norm preserving}, we have
	    \begin{align}
	      \big\|S_K^* f^{\lambda}_{\b{z}} - f_{\rho}\big\|_{\rho}  &=  \big\|S_K^*\big(f^{\lambda}_{\b{z}} - f^{\lambda}\big) +       S_K^*f^{\lambda} - f_{\rho}\big\|_{\rho}
	      \le  \big\|S_K^*\big(f^{\lambda}_{\b{z}} - f^{\lambda}\big)\big\|_{\rho} +       \big\|S_K^*f^{\lambda} - f_{\rho}\big\|_{\rho}\nonumber \\
	      &= \big\|\sqrt{T}\big(f^{\lambda}_{\b{z}} - f^{\lambda}\big)\big\|_{\H} +       \big\|S_K^*f^{\lambda} - f_{\rho}\big\|_{\rho} \label{eq:fzlambda-frho:bound}.
	    \end{align}
      \item  \tb{Decomposition of $\big\|\sqrt{T}\left(f^{\lambda}_{\b{z}} - f^{\lambda}\right)\big\|_{\H}$}:
	    Making use of the analytical expressions for $f^{\lambda}_{\b{z}}$ and $f^{\lambda}$ [see Eq.~\eqref{eq:f_zlambda} and Eq.~\eqref{eq:flambda:def}],
	    and the operator Woodbury formula \citep[][Theorem~2.1, page 724]{ding08spectrum} we arrive at the decomposition (see Section~\ref{proof:fzlambda-flambda:decomposition})
	    \begin{align*}
		\big\|\sqrt{T}\big(f^{\lambda}_{\b{z}} - f^{\lambda}\big)\big\|_{\H} 
		&\le \big\|\sqrt{T}(T_{\b{x}}+\lambda I)^{-1}\big\|_{\L(\H)}  \Big( \left\|g_{\b{z}}  - g_{\rho}\right\|_{\H} + \vphantom{\left\| T-T_{\b{x}}\right\|_{\L(\H)} \lambda^{-1} \left\| S_K\left[f_{\rho} - \left(\tilde{T}+\lambda I\right)^{-1}S_K^*S_Kf_{\rho} \right] \right\|_{\H}}\Big.\nonumber\\
		& \hspace*{1.3cm} \Big. \left\| T-T_{\b{x}}\right\|_{\L(\H)} \lambda^{-1} \big\| S_K\big[f_{\rho} - (\tilde{T}+\lambda I)^{-1}S_K^*S_Kf_{\rho} \big] \big\|_{\H} \Big),
	    \end{align*}
	    where $g_{\rho} = S_Kf_{\rho}$. As it is known \citep[page~348]{caponnetto07optimal} $\| \sqrt{T} (T_{\b{x}}+\lambda I)^{-1}\|_{\L(\H)} \le 1/\sqrt{\lambda}$
	    provided that  $\bm{\Theta}(\lambda,\b{z}) \le \frac{1}{2}$.
    \item \tb{Bound on $\left\|g_{\b{z}}  - g_{\rho}\right\|_{\H}$, $\left\|T-T_{\b{x}}\right\|_{\L(\H)}$}: By concentration arguments the bounds
	  \begin{align*}
	      \left\|g_{\b{z}}  - g_{\rho}\right\|_{\H} \le \left( \frac{4C\sqrt{B}_K}{l}+ \frac{2C\sqrt{B_K}}{\sqrt{l}} \right)\log\left(\frac{2}{\eta}\right),\hspace*{0.1cm}
	      \left\| T-T_{\b{x}}\right\|_{\L(\H)} \le \left(\frac{4B_K}{l}+\frac{4\sigma}{\sqrt{l}}\right)\log\left(\frac{2}{\eta}\right)
	  \end{align*}
	  hold with probability at least $1-\eta$, each (see Section~\ref{proof:gz-grho:bound}, \ref{proof:T-Tx:bound}).
    \item \tb{Decomposition of $\big\| S_K\big[f_{\rho} - (\tilde{T}+\lambda I)^{-1}S_K^*S_Kf_{\rho} \big] \big\|_{\H}^2$}: Exploiting the analytical
	  formula for $f^{\lambda}$, one can construct (Section~\ref{proof:S_Klong}) the upper bound
	  \begin{align*}
			\big\| S_K\big[f_{\rho} - (\tilde{T}+\lambda I)^{-1}S_K^*S_Kf_{\rho} \big] \big\|_{\H}^2 &\le \big\| \tilde{T}\big[f_{\rho} - (\tilde{T}+\lambda I)^{-1}S_K^*S_Kf_{\rho} \big] \big\|_{\rho}
			\big\|S_K^*f^{\lambda} - f_{\rho}\big\|_{\rho}.
	  \end{align*}
    \item \tb{Bound on $\big\| \tilde{T}\big[f_{\rho} - (\tilde{T}+\lambda I)^{-1}S_K^*S_Kf_{\rho} \big] \big\|_{\rho}$}: 
	    Using our  assumptions that $f_{\rho}\in \Im(\tilde{T}^s)$ ($s\ge 0$)\footnote{Note that we choose $s=0$ and $s>0$ in the first and second theorem part, respectively.} and exploiting the separability of $L^{2}_{\rho_X}$, by 
	    Lemma~\ref{lemma:lin trafo-ed:ridge} ($\K = L^2_{\rho_X}$, $f=f_{\rho}$, $M=\tilde{T}$, $a=1$) and $\tilde{T}=S_K^*S_K$  we obtain the upper bound 
			\begin{eqnarray*}
			    \lefteqn{\big\| \tilde{T}\big[f_{\rho} - (\tilde{T}+\lambda I)^{-1}S_K^*S_Kf_{\rho} \big] \big\|_{\rho} = \big\| \tilde{T}\big[f_{\rho} - (\tilde{T}+\lambda I)^{-1}\tilde{T}f_{\rho} \big] \big\|_{\rho}}\\
				& \le   \max\Big(1,\|\tilde{T}\|^s_{\L(L^2_{\rho_X})}\Big) \lambda^{\min(1,s+1)}\big\|\tilde{T}^{-s}f_{\rho}\big\|_{\rho} =   \max\Big(1,\|\tilde{T}\|^s_{\L\left(L^2_{\rho_X}\right)}\Big) \lambda\big\|\tilde{T}^{-s}f_{\rho}\big\|_{\rho},
			\end{eqnarray*}
			where we used at the last step that $\min(1,s+1)=1$; this follows from $s\ge 0$.
    \item \tb{Bound on $\left\|S_K^*f^{\lambda} - f_{\rho}\right\|_{\rho}$}:
	  \begin{enumerate}
	    \item \tb{No range space assumption}: One can construct (Section~\ref{proof:flambda-frho-with-h}) the bound 
			\begin{align*}
			    \big\|S_K^*f^{\lambda} - f_{\rho}\big\|_{\rho} &\le  \left\|f_{\rho}-S_K^*q\right\|_{\rho} + \max\big(1,\left\|T\right\|_{\L(\H)}\big)\lambda^{\frac{1}{2}} \left\|q\right\|_{\H},
			\end{align*}
			which holds for arbitrary $q\in \H$.
	    \item \tb{Range space assumption in $L^2_{\rho_X}$}: Using the $S_K^* f^{\lambda} = (\tilde{T}+\lambda I)^{-1}\tilde{T} f_{\rho}$ identity [see Eq.~\eqref{eq:SKf-lambda:conversion}], and Lemma~\ref{lemma:lin trafo-ed:ridge} ($M=\tilde{T}$, $\K=L^2_{\rho_X}$, $a=0$), we get
		\begin{align*}
		    \big\|S_K^*f^{\lambda} - f_{\rho}\big\|_{\rho} &= \big\| (\tilde{T}+\lambda I)^{-1}\tilde{T} f_{\rho} - f_{\rho}\big\|_{\rho}
		    \le \max\Big(1,\big\|\tilde{T}\big\|_{\L(L^2_{\rho_X})}^{s-1}\Big) \lambda^{\min(1,s)}\big\|\tilde{T}^{-s}f_{\rho}\big\|_{\rho}.
		\end{align*}
	  \end{enumerate}
    \item \tb{Union bound}: Applying an $\alpha=\log(l)+\delta$ reparameterization, changing $\eta$ to $\frac{\eta}{3}$ and combining the derived results (in case of the first statement with $s=0$) with a union bound, Theorem~\ref{theo:L2} follows.
  \end{enumerate}

\begin{remark}
\label{remark:well-vs-misspec} To contrast the derivation of the well- and the misspecified cases, we note that previous results [Section~\ref{sec:results:A1}, or \citet{caponnetto07optimal}'s bound] were used at two points:
	\begin{itemize}
	  \item[(a)] In Step~2 by using Eq.~\eqref{eq:S_K^*:norm preserving} and transforming the $L^2_{\rho_X}$ error $\left\|S_K^* \left(f^{\lambda}_{\hat{\b{z}}} - f^{\lambda}_{\b{z}}\right)\right\|_{\rho}$ to
	      $\H$, we could rely on our previous bounds for  $S_{-1}$ and $S_0$. However, we were required to use a different concentration argument to guarantee $\bm{\Theta}(\lambda,\b{z}) \le \frac{1}{2}$ since  we no longer assume the 
	      $\P(b,c)$ prior class.
	  \item[(b)] In Step~4 the first term could be bounded by \citet{caponnetto07optimal}. Its $\bm{\Theta}(\lambda,\b{z}) \le \frac{1}{2}$ condition was  guaranteed by Step~2; and see Section~\ref{proof:S_-1+S_0}.
	  \end{itemize}
      We note that our misspecified proof method was inspired by \citet[Theorem 12]{sriperumbudur14density}, where the authors focused on the consistency of an infinite-dimensional exponential family estimator.
\end{remark}

\section{Related Work} \label{sec:related-work}
In this section we discuss existing approaches and heuristic techniques to tackle learning problems on distributions.

\vspace{1mm}\tb{Methods based on parametric assumptions:} A number of methods have been proposed to compute the similarity of distributions or bags of samples. As a first approach, one could fit
a parametric model to the bags, and estimate the similarity of the bags based on the obtained parameters. It is then possible to
define learning algorithms on the basis of these similarities, which often take analytical form. Typical examples with explicit formulas include Gaussians,
finite mixtures of Gaussians, and distributions from the exponential family \citep[with known log-normalizer function and zero carrier measure, see][]{kondor03kernel,jebara04probability,wang09closed,nielsen12closed}.
A major limitation of these methods, however, is that they apply quite simple parametric assumptions, which may not be sufficient or verifiable in practise.\vspace{1mm}

\tb{Methods based on parametric assumption in a RKHS:} A heuristic related to the parametric approach is to assume that the training distributions are Gaussians in a reproducing kernel Hilbert space \citep[see for
example][and references therein]{jebara04probability,zhou06sample}. This assumption is algorithmically appealing, as many divergence measures for Gaussians can be
computed in closed form using only inner products, making them straightforward to kernelize. A fundamental shortfall of kernelized Gaussian divergences is the lack of
their consistency analysis in specific learning algorithms.\vspace{1mm}

\tb{Kernels based techniques:} A more theoretically grounded approach to learning on distributions has been to define positive definite kernels on the basis of statistical divergence measures on distributions,
or by metrics on non-negative numbers; these can then be used in kernel algorithms.
This category includes work on semigroup kernels \citep{cuturi05semigroup}, non-extensive information theoretical kernel constructions \citep{martins09nonextensive}, and kernels based
on Hilbertian metrics \citep{hein05hilbertian}. For example, the intuition of semigroup kernels \citep{cuturi05semigroup} is as follows: if two measures or sets of points
overlap, then their sum is expected to be more concentrated. The value of dispersion can be measured by entropy or inverse generalized variance. In the second type of approach
\citep{hein05hilbertian}, homogeneous Hilbert metrics on the non-negative real line are used to define the similarity of probability distributions.
While these techniques guarantee to provide valid kernels on certain restricted domains of measures, the performance of learning
algorithms based on finite-sample estimates of these kernels remains a challenging open question.
One might also plug into learning algorithms (based on similarities of distributions) consistent R{\'e}nyi and Tsallis divergence estimates \citep{poczos11nonparametric,poczos12support}, but
these similarity indices are \emph{not} kernels, and their consistency in specific \emph{learning tasks}
remains an open question.
\vspace{1mm}

\tb{Multi-instance learning:} An alternative paradigm in learning when the inputs are ``bags of objects'' is to simply treat each input as a \emph{finite
set}: this leads to the multi-instance learning task \citep[MIL, see][]{dietterich97solving,ray01multiple,dooly02multiple}. 
In MIL one is given a set of labelled bags, and the task of the learner is to
find the mapping from the bags to the labels. Many important examples fit into the MIL framework: for example, different configurations of a given molecule can 
be handled as a bag of shapes, images can be considered as a set of patches or regions of interest, a video can be seen as a collection of images, 
a document might be described as a bag of words or paragraphs, a web page can be identified by its links, a group of people on a social network can be captured by their friendship graphs, 
in a biological experiment a subject can be identified by his/her time series trials, or a customer might be characterized by his/her shopping records.
The MIL approach has been applied in
several domains; see the reviews from \citet{babenko04multiple,zhou04multi,foulds10review,amores13multiple}. \vspace{1mm}

\tb{``Bag-of-objects'' methods (MIL, classification):} Despite the large number of MIL applications and the spate of heuristic solution techniques, there exist few \emph{theoretical results} in the area \citep{auer98approximating,long98pac,blum98note,babenko11multiple,zang13milage,sabato12multi} and they focus on the
multi-instance \emph{classification} (MIC) task. In particular, let us first consider the standard MIC assumption \citep{dietterich97solving}, where a bag is declared to be positive (labelled with ``1'') if at least one of its instances is positive (``1''); 
otherwise, the bag is negative (``0'').\footnote{The motivation of this assumption comes from drug discovery: if a molecule has at least one well-binding configuration, then it is considered to bind well.}
In other words, if the instances ($x_{i,n}$) in the $i^{th}$ bag $\{x_{i,1},\ldots,x_{i,N}\}$ have hidden label $L(x_{i,n})\in\{0,1\}$, then the observed label of the bag is 
$y_i=h(x_{i,1},\ldots,x_{i,N})=\max(L(x_{i,1}),\ldots,L(x_{i,N}))\in\{0,1\}$. In case of the original APR \citep[axis-aligned rectangles;][]{dietterich97solving} hypothesis class,  
function $L$ is equal to the indicator of an unknown rectangle $R$ ($L=\I_R$). 
In other words, a bag is declared to be positive if there exists at least one instance in the bag, which belongs to $R$.\footnote{In terms of drug binding prediction, this means that a molecule binds to a target iff at least one of its configurations falls within 
a fixed, but unknown rectangle.} The goal is to learn $R$ with high probability given the bags ($\{x_{i,1},\ldots,x_{i,N}\}$-s) and their labels ($y_i$-s).
\citet{long98pac} proved the PAC learnability \citep[probably approximately correct;][]{valiant84theory} of the APR hypothesis class, if 
all instances in each bag are i.i.d.\ and follow the same product distribution over the instance coordinates. On the other hand, for arbitrary distributions over bags, when 
the instances within a bag might be statistically dependent, APR learning under MIC is NP-hard 
\citep{auer98approximating}; the same authors also showed that the product property \citep{long98pac} on the coordinates is not required to obtain PAC results. 
\citet{blum98note} extended PAC learnability of APR-s to hypothesis classes learnable from one-sided classification noise.
In contrast to the previous approaches \citep{long98pac,auer98approximating,blum98note}, \citet{babenko11multiple} modelled the bags as 
low-dimensional manifolds, and proved PAC bounds. By relaxing the standard MIC assumption, \citet{sabato12multi} showed PAC-learnability for general MIC hypothesis classes with extended
 ``$\max$'' functions. \citet{zang13milage} derived high-probability generalization bounds in the MIC setting, when local and global representations are combined.
Our work falls outside this setting 
since  the label and bag generation mechanisms we consider are different: we do not assume an exact form of the labelling mechanism (function $L$ and $\max$ in $h$). Rather, the labelling is presumed to be stochastically determined by the underlying true distribution, not deterministically by the instance realizations in the bags (these are presumed i.i.d., and may be bag-specific).
\vspace{1mm}

\tb{``Bag-of-objects'' methods (MIL, not classification):} Beyond classification, there exist several \emph{heuristics}---without consistency guarantees---for many other multi-instance problems in the literature, including
 regression \citep{ray01multiple,dooly02multiple,zhou09multi,kwok07marginalized}, clustering \citep{zhang09multi,zhang09m3ic,zhang11maximum,chen12contextual}, 
 ranking \citep{bergeron08multiple,hu08multiple,bergeron12fast}, 
 outlier detection \citep{wu10identifying}, transfer learning \citep{raykar08bayesian,zhang09multiple}, and feature selection, -weighting and -extraction 
\citep[also called dimensionality reduction, low-dimensional embedding, manifold learning, see][and references therein]{raykar08bayesian,ping10noniid,sun10multi,carter11information,zafrahybrid13,chai14maximum,chai14multiple}.\vspace{1mm}

\tb{Approaches using set metrics:} Adapting the bag viewpoint of MIL, one can come up with set metric based learning algorithms.\footnote{Often these ``metrics'' are only semi-metrics, as they do not satisfy the triangle inequality.}
Probably one of the most well-known set metrics is the Hausdorff metric \citep{edgar95measure}, which is defined for non-empty compact sets of metric spaces, specifically for sets
containing finitely many points. There also exist other (semi)metric constructions on points sets \citep{eiter97distance,ramon01polynomial}.
Unfortunately, the classical Hausdorff metric is highly sensitive to outliers, seriously limiting its practical applicability. In order to mitigate this deficiency, several variants of the 
Hausdorff metric  have been designed in the MIL literature, such as  the maximal-, the minimal- and the ranked Hausdorff metrics, with successful applications in 
MIC \citep{wang08solving} and multi-instance outlier detection \citep{wu10identifying};
and the average Hausdorff metric \citep{zhang09multi} and  contextual Hausdorff dissimilarity  \citep{chen12contextual}, which have been found useful in multi-instance clustering. 
Unfortunately,  these methods lack  any theoretical guarantee when applied in specific learning problems.\vspace{1mm}

\tb{Functional data analysis techniques:} Finally, the distribution regression task might also be interpreted as a functional data analysis problem \citep{ramsay02applied,ramsay05functional,muller05functional}, by considering the probability measures $x_i$  as 
functions. This is a highly non-standard setup, however, since these functions ($x_i$) are defined on $\sigma$-algebras and are non-negative, $\sigma$-additive.

\section{Conclusion}\label{sec:conclusions}
We have established a learning theory of distribution regression,  where the inputs are
probability measures on separable, topological domains endowed with reproducing kernels,
and the outputs are elements of a separable Hilbert space.
We studied a ridge regression scheme defined on embeddings of the input distributions to a reproducing kernel Hilbert space,
which has a simple analytical solution, as well as theoretically sound, efficient methods for 
approximation \citep{zhang15divide,richtarik16distributed,alaoui15fast,yang16randomized,rudi15less}. We derived explicit bounds on the excess risk as a function of the number of samples  and problem difficulty. We tackled both the well-specified case (when the regression function belongs to the assumed RKHS modelling class), 
and the more general misspecified setup. As a special case of our results, we proved the consistency of regression for set kernels \citep{haussler99convolution,gartner02multi}, which was a $17$-year-old open problem, 
and  for a recent kernel family \citep{christmann10universal}, which we have expanded upon (Table \ref{tab:K examples}). 
We proved an exact computational-statistical efficiency trade-off for the MERR estimator: in the 
well-specified setting, we showed how to choose the bag size in the two-stage sampled setup to match the one-stage sampled minimax optimal rate \citep{caponnetto07optimal}; and in the misspecified setting, our rates
approximate closely an asymptotically optimal estimator  imposing stricter eigenvalue decay conditions \citep{steinwart09optimal}.
Several exciting open questions remain, including whether  improved/optimal rates can be derived in the misspecified case, whether we can obtain consistency guarantees for non-point estimates, and how to handle non-ridge extensions. 

Finally, we note that although the primary focus of the current paper was theoretical, we have applied the 
MERR method \citep[Section~A.2]{szabo15twostage} to
supervised entropy learning and aerosol prediction based on multispectral satellite images.\footnote{For code, see \url{https://bitbucket.org/szzoli/ite/}.}  In future work, we will address applications with vector-valued outputs.

\section{Proofs}\label{sec:proofs}
We provide proofs for our results detailed in Section~\ref{sec:convergence analysis}: Section~\ref{sec:proofs:well-specified} (\emph{resp.} Section~\ref{sec:proofs:misspecified}) focuses on the well-specified case (\emph{resp.} misspecified setting). 
The used lemmas are enlisted in Section~\ref{sec:supplement}.

\subsection{Proofs of the Well-specified Case}\label{sec:proofs:well-specified}
We give proof details concerning the excess risk in the well-specified case (Theorem~\ref{theo1}).

\subsubsection{Proof of the bound on $\|  g_{\hat{\b{z}}} - g_{\b{z}} \|_{\H}^2$}\label{proof:gz}
	    By \eqref{eq:f_zlambda}, \eqref{eq:f_hatz^lambda} we get $g_{\hat{\b{z}}} - g_{\b{z}} = \frac{1}{l}\sum_{i=1}^l \big( K_{\mu_{\hat{x}_i}} - K_{\mu_{x_i}} \big) y_i$; hence by applying the H{\"o}lder property of $K_{(\cdot)}$, the boundedness of 
	    $y_i$ ($\left\|y_i\right\|_Y\le C$)
	    and  \eqref{eq:emp-mean-emb-conv-rate:Assumption1}, we obtain
	    \begin{eqnarray*}
		    \lefteqn{\left\|  g_{\hat{\b{z}}} - g_{\b{z}}  \right\|_{\H}^2  \le  \frac{1}{l^2} l \sum_{i=1}^l \big\|  \big(K_{\mu_{\hat{x}_i}} - K_{\mu_{x_i}}\big) y_i \big\|_{\H}^2
		      \le  \frac{1}{l} \sum_{i=1}^l \big\| K_{\mu_{\hat{x}_i}} - K_{\mu_{x_i}} \big\|_{L(Y,\H)}^2 \left\|y_i\right\|_{Y}^2}\nonumber\\
		      &&\hspace*{-0.5cm}\le \frac{L^2}{l} \sum_{i=1}^l \left\|y_i\right\|_Y^2 \left\|\mu_{\hat{x}_i}-\mu_{x_i}\right\|_{H}^{2h}
		      \le \frac{L^2C^2}{l}\sum_{i=1}^l\left[ \frac{\left(1+\sqrt{\alpha}\right)\sqrt{2B_k}}{\sqrt{N}} \right]^{2h}
		      =L^2 C^2  \frac{\left(1+\sqrt{\alpha}\right)^{2h} (2B_k)^h}{N^h} 
	    \end{eqnarray*}
	    with probability at least $1 - l e^{-\alpha}$, based on a union bound.

\subsubsection{Proof of the bound on $\| T_{\b{x}}-T_{\hat{\b{x}}}\|_{\L(\H)}^2$}\label{proof:Tx}
	 Using the definition of $T_{\b{x}}$ and $T_{\hat{\b{x}}}$, and  exploiting (with $\|\cdot\|_{\L(\H)}$) that in a  
	 normed space\footnote{Eq.~\eqref{eq:norm-eq} holds since $\left\|\cdot\right\|^2$ is 
	 convex function, thus $\left\|\frac{1}{n}\sum_{i=1}^nf_i\right\|^2\le \frac{1}{n}\sum_{i=1}^n\left\|f_i\right\|^2$.\label{footnote:norm-eq-reasoning}} $(N,\left\|\cdot\right\|)$, $f_i\in N$, ($i=1,\ldots, n$)
	\begin{align}
		    \big\|\sum\nolimits_{i=1}^nf_i\big\|^2 & \le n \sum\nolimits_{i=1}^n\left\|f_i\right\|^2,  \label{eq:norm-eq}
	 \end{align}
          we get
	  \begin{align}
	      \left\|T_{\b{x}} - T_{\hat{\b{x}}}\right\|_{\L(\H)}^2 &\le \frac{1}{l^2} l \sum_{i=1}^l \left\| T_{\mu_{x_i}} - T_{\mu_{\hat{x}_i}}\right\|_{\L(\H)}^2. \label{eq:Tx-Txi}
	  \end{align}
	  To upper bound $\| T_{\mu_{x_i}} - T_{\mu_{\hat{x}_i}}\|_{\L(\H)}^2$, let us see how $T_{\mu_u}=K_{\mu_a}K^*_{\mu_a}$ acts.
	  The existence of an $E\ge 0$  constant satisfying $\left\|(T_{\mu_u}-T_{\mu_v})(f)\right\|_{\H} \le E \left\|f\right\|_{\H}$
	    implies $\left\|T_{\mu_u}-T_{\mu_v}\right\|_{\L(\H)} \le E$. We continue with the l.h.s.\ of this equation using Eq.~\eqref{eq:norm-eq}: 
		  \begin{align*}
			\left\|(T_{\mu_u}-T_{\mu_v})(f)\right\|_{\H}^2
			  &=\left\| K_{\mu_u} K_{\mu_u}^*(f) - K_{\mu_v} K_{\mu_v}^*(f)\right\|_{\H}^2\\
			  &= \left\| K_{\mu_u } \left[K_{\mu_u}^*(f) -  K_{\mu_v}^*(f)\right]
			  +\left( K_{\mu_u} - K_{\mu_v}\right) K_{\mu_v}^*(f) \right\|_{\H}^2 \nonumber\\
			  &\le 2 \left[\left\| K_{\mu_u } \left[K_{\mu_u}^*(f) -  K_{\mu_v}^*(f)\right] \right\|_{\H}^2
			  + \left\|\left( K_{\mu_u} - K_{\mu_v}\right) K_{\mu_v}^*(f) \right\|_{\H}^2\right]. \nonumber
		  \end{align*}
		    By  Eq.~\eqref{eq:phiK:bounded} and the H{\"o}lder continuity of $K_{(\cdot)}$, one arrives at
		  \begin{eqnarray*}
		      \lefteqn{\left\| K_{\mu_u } \left[K_{\mu_u}^*(f) -  K_{\mu_v}^*(f)\right] \right\|_{\H}^2  \le \left\|K_{\mu_u }\right\|_{\L(Y,\H)}^2 \left\|K_{\mu_u}^*(f) -  K_{\mu_v}^*(f)\right\|_{Y}^2}\\
			    &&\le \left\|K_{\mu_u }\right\|_{\L(Y,\H)}^2 \left\|K_{\mu_u}^* -  K_{\mu_v}^*\right\|_{\L(\H,Y)}^2 \left\|f\right\|_{\H}^2
			    = \left\|K_{\mu_u }\right\|_{\L(Y,\H)}^2 \left\|\left(K_{\mu_u} -  K_{\mu_v}\right)^*\right\|_{\L(\H,Y)}^2 \left\|f\right\|_{\H}^2\\
			    &&= \left\|K_{\mu_u }\right\|_{\L(Y,\H)}^2 \left\|K_{\mu_u} -  K_{\mu_v}\right\|_{\L(Y,\H)}^2 \left\|f\right\|_{\H}^2
			    \le B_K L^2\left\|\mu_u -  \mu_v\right\|_{H}^{2h} \left\|f\right\|_{\H}^2,\\
			  \lefteqn{\left\|\left( K_{\mu_u} - K_{\mu_v}\right) K_{\mu_v}^*(f) \right\|_{\H}^2  \le \left\| K_{\mu_u} - K_{\mu_v} \right\|_{\L(Y,\H)}^2 \left\|K_{\mu_v}^*(f)\right\|_{Y}^2}\\ 
			    &&\le \left\| K_{\mu_u} - K_{\mu_v} \right\|_{\L(Y,\H)}^2 \left\|K_{\mu_v}^*\right\|_{\L(\H,Y)}^2 \left\|f\right\|_{\H}^2
			    \le B_K L^2\left\|\mu_u -  \mu_v\right\|_{H}^{2h} \left\|f\right\|_{\H}^2.
		  \end{eqnarray*}
		  Hence $\left\|(T_{\mu_u}-T_{\mu_v})(f)\right\|_{\H}^2 \le 4 B_K L^2\left\|\mu_u -  \mu_v\right\|_{H}^{2h} \left\|f\right\|_{\H}^2\hspace{0.1cm} \Rightarrow\hspace{0.1cm}
		    E^2 = 4 B_K L^2\left\|\mu_u -  \mu_v\right\|_{H}^{2h}$. Exploiting this property in \eqref{eq:Tx-Txi}  with Eq.~\eqref{eq:emp-mean-emb-conv-rate:Assumption1} we arrive to the bound
		\begin{align}
		      \left\|T_{\b{x}} - T_{\hat{\b{x}}}\right\|_{\L(\H)}^2
		      & \le  \frac{4B_KL^2}{l}\sum_{i=1}^l\left\|\mu_{x_i}-\mu_{\hat{x}_i}\right\|_H^{2h}
		      \le \frac{4B_K L^2}{l}\sum_{i=1}^l \frac{\left(1+\sqrt{\alpha}\right)^{2h} (2B_k)^h}{N^h} \nonumber\\
		      &= \frac{\left(1+\sqrt{\alpha}\right)^{2h} 2^{h+2}(B_k)^{h}B_K L^2}{N^h}.\label{eq:Tx-Txhat}
		\end{align}

\subsubsection{Proof: final union bound in Theorem~\ref{theo1}}\label{proof:theo1:union bound}
    Until now, we obtained that if (i) the sample number $N$ satisfies Eq.~\eqref{eq:N}, (ii)
	  \eqref{eq:emp-mean-emb-conv-rate:Assumption1} holds (which has probability at least $1 - l e^{-\alpha} = 1 - e^{-[\alpha-\log(l)]}= 1-e^{-\delta}$ applying a union bound argument; $\alpha=\log(l)+\delta$), and
           (iii) $\bm{\Theta}(\lambda,\b{z})\le \frac{1}{2}$ is fulfilled [see Eq.~\eqref{eq:theta}], then
  \begin{eqnarray*}
	\lefteqn{S_{-1} + S_{0}  \le \frac{4}{\lambda} \left[L^2 C^2  \frac{\left(1+\sqrt{\alpha}\right)^{2h} (2B_k)^h}{N^h}
	       + \frac{\left(1+\sqrt{\alpha}\right)^{2h} 2^{h+2}(B_k)^{h}B_K L^2}{N^h}  \times \right.}\\
		 &&\left.\hspace*{-0.7cm} \times \left(\log^2\left(\frac{6}{\eta}\right) \left\{\frac{64}{\lambda} \left[\frac{M^2 B_K}{l^2\lambda} + \frac{\Sigma^2\N(\lambda)}{l} \right]
		   + \frac{24}{\lambda^2} \left[ \frac{4B_K^2\B(\lambda)}{l^2} + \frac{B_K \A(\lambda)}{l}\right]\right\}  + \B(\lambda) + \left\| f_{\rho}\right\|_{\H}^2\right)\right] \nonumber\\
	      &&\hspace*{0.95cm}= \frac{4L^2\left(1+\sqrt{\alpha}\right)^{2h}(2B_k)^h}{\lambda N^h}\left[C^2+4B_K\times \phantom{ \left(\log^2\left(\frac{6}{\eta}\right) \left\{\frac{64}{\lambda} \left[\frac{M^2 B_K}{l^2\lambda} + \frac{\Sigma^2\N(\lambda)}{l} \right]
		   + \frac{24}{\lambda^2} \left[ \frac{4B_K^2\B(\lambda)}{l^2} + \frac{B_K \A(\lambda)}{l}\right]\right\}  + \B(\lambda) + \left\| f_{\rho}\right\|_{\H}^2\right)}\right.\\
		&& \left.\hspace*{-0.7cm} \times \left(\log^2\left(\frac{6}{\eta}\right) \left\{\frac{64}{\lambda} \left[\frac{M^2 B_K}{l^2\lambda} + \frac{\Sigma^2\N(\lambda)}{l} \right]
		   + \frac{24}{\lambda^2} \left[ \frac{4B_K^2\B(\lambda)}{l^2} + \frac{B_K \A(\lambda)}{l}\right]\right\}  + \B(\lambda) + \left\| f_{\rho}\right\|_{\H}^2\right)\right].\nonumber
  \end{eqnarray*}
  By taking into account \citet{caponnetto07optimal}'s bounds for $S_1$ and $S_2$, 
    $S_1 \le 32 \log^2\left(\frac{6}{\eta}\right)\left[\frac{B_K M^2}{l^2\lambda} + \frac{\Sigma^2 \N(\lambda)}{l} \right]$, 
    $S_2 \le 8 \log^2\left(\frac{6}{\eta}\right) \left[\frac{4B_K^2\B(\lambda)}{l^2\lambda} + \frac{B_K \A(\lambda)}{l\lambda} \right]$,
  plugging all the expressions to \eqref{eq:E:5terms}, we obtain Theorem~\ref{theo1} with a union bound.

\subsection{Proofs of the Misspecified Case}\label{sec:proofs:misspecified}
We present the proof details concerning the excess risk in the misspecified case (Theorem~\ref{theo:L2}).

\subsubsection{Proof of the bound on $\sqrt{S_{-1}} + \sqrt{S_0}$ without $\P(b,c)$} \label{proof:S_-1+S_0}
	The upper bounds on $S_{-1}$ and $S_0$ [which are defined in Eqs.~\eqref{eq:Assumption1:S_-1:def}, \eqref{eq:Assumption1:S_0:def}] remain valid without modification
	provided that
        (i) $\bm{\Theta}(\lambda,\b{z}) = \|(T - T_{\b{x}}) (T+\lambda)^{-1}\|_{\L(\H)} \le \|(T - T_{\b{x}}) (T+\lambda)^{-1}\|_{\L_2(\H)} \le \frac{1}{2}$, where we used Eq.~\eqref{eq:opnorm<=HSnorm},
        (ii) Eq.~\eqref{eq:emp-mean-emb-conv-rate:Assumption1} is satisfied (which has probability $1-le^{-\alpha}$) and (iii) Eq.~\eqref{eq:N} holds.
	Our goal below is to guarantee the $\bm{\Theta}(\lambda,\b{z}) \le \frac{1}{2}$ condition with high probability \emph{without} assuming that the prior belongs to $\P(b,c)$.

	\noindent\tb{Requirement $\bm{\Theta}(\lambda,\b{z}) \le \frac{1}{2}$}: Let us define $\xi_i = T_{\mu_{\b{x}_i}}(T+\lambda)^{-1}\in \L_2(\H)$, $(i=1,\ldots,l)$. With this choice we get $\E[\xi_i] = T(T+\lambda)^{-1}$, $(T - T_{\b{x}})(T+\lambda)^{-1} = \E[\xi_i] - \frac{1}{l}\sum_{i=1}^l\xi_i$ and 
	\begin{align*}
	  \left\|\xi_i\right\|_{\L_2(\H)} &\le \left\|T_{\mu_{\b{x}_i}}\right\|_{\L_2(\H)}  \left\|(T+\lambda)^{-1}\right\|_{\L(\H)} \le B_K/\lambda \hspace{0.1cm} \Rightarrow \hspace{0.1cm}
	  \E\big[\left\|\xi_i\right\|_{\L_2(\H)}^2\big] \le (B_K)^2 / \lambda^2, %\label{eq:xi-norm}
	\end{align*}
	where we made use of \eqref{eq:prod-norm-ineq}, the $\|T_{\mu_{\b{x}_i}}\|_{\L_2(\H)}\le B_K$ identity following from the boundedness of $K$ \citep[][page 341, Eq.~(13)]{caponnetto07optimal}, and the spectral theorem. Consequently, by
	the Bernstein's inequality (Lemma~\ref{lemma:Bernstein} with $\K=\L_2(\H)$, $B=2B_K/\lambda$, $\sigma=B_K/\lambda$) we obtain that for $\forall \eta\in (0,1)$
	\begin{align*}
	  \Pr\left(\left\|(T - T_{\b{x}})(T+\lambda)^{-1}\right\|_{\L_2(\H)} \le 2\left(\frac{2B_K}{\lambda l}+\frac{B_K}{\sqrt{l}\lambda}\right) \log\left(\frac{2}{\eta}\right)\right) \ge 1-\eta.
	\end{align*}
	Thus, for $\bm{\Theta}(\lambda,\b{z}) \le \frac{1}{2}$ with probability $1-\eta$ it is sufficient to have
	\begin{align}
	  2\left(\frac{2B_K}{\lambda l}+\frac{B_K}{\sqrt{l}\lambda}\right) \log\left(\frac{2}{\eta}\right) & \le \frac{6 B_K}{\sqrt{l} \lambda} \log\left(\frac{2}{\eta}\right) \le \frac{1}{2} 
	  \Leftrightarrow \left[\frac{12 B_K}{\lambda} \log\left(\frac{2}{\eta}\right)\right]^2 \le l. \label{eq:r2}
	\end{align}
	
	Under these conditions, we arrived at the upper bound
	\begin{align*}
	  \sqrt{S_{-1}} + \sqrt{S_0} &\le \sqrt{\frac{4L^2 C^2 \left(1+\sqrt{\alpha}\right)^{2h} (2B_k)^h}{\lambda N^h}} \left[ \sqrt{1} + \sqrt{\frac{4 B_K}{\lambda}}\right]\\
	  &= \frac{2LC(1+\sqrt{\alpha})^h(2B_k)^{\frac{h}{2}}}{\sqrt{\lambda}N^{\frac{h}{2}}} \left[1 + \frac{2\sqrt{B_K}}{\sqrt{\lambda}}\right],
	\end{align*}
	where as opposed to Section~\ref{proof:theo1:union bound} and Eq.~\eqref{f:zlambda:improved-bound} we used a slightly cruder $\left\|f_{\b{z}}^{\lambda}\right\|_{\H}^2\le \frac{C^2}{\lambda}$ bound; 
	  it holds  without the $\P(b,c)$ assumption by the definition of $f_{\b{z}}^{\lambda}$ and the boundedness of $y$ since $\lambda\left\|f_{\b{z}}^{\lambda}\right\|_{\H}^2\le \frac{1}{l}\sum_{i=1}^l \left\|y_i\right\|_Y^2 \le C^2$.\vspace{1mm}

        \tb{Remark}:\label{remark:without-N(lambda)-price} Notice that the price we pay for not assuming that the prior belongs to the $\P(b,c)$ class ($b>1$) is a slightly tighter 
	$\frac{1}{\lambda^2}\le l$ constraint [Eq.~\eqref{eq:r2}] instead of $\frac{1}{\lambda^{1+\frac{1}{b}}}\le l$ in Eq.~\eqref{eq:f1}, and a somewhat looser $\left\|f_{\b{z}}^{\lambda}\right\|_{\H}^2$ bound.

\subsubsection{Proof of the decomposition of $\big\|\sqrt{T}\left(f^{\lambda}_{\b{z}} - f^{\lambda}\right)\big\|_{\H}$}\label{proof:fzlambda-flambda:decomposition}
	    Using the analytical formula of $f^{\lambda}_{\b{z}}$ [see Eq.~\eqref{eq:f_zlambda}] and that of $f^{\lambda}$ [see Eq.\eqref{eq:flambda:def}]
	      \begin{align}
		f^{\lambda}&=(S_KS_K^*+\lambda I)^{-1}S_Kf_{\rho}=(T+\lambda I)^{-1}S_Kf_{\rho}  \label{eq:analytical expr:f-lambda}     
	      \end{align}
	    one gets
	    $(T+\lambda I)f^{\lambda}=S_Kf_{\rho}$ $\Rightarrow$ $\lambda f^{\lambda} = S_Kf_{\rho} - Tf^{\lambda}$ and
	    \begin{align}
		f^{\lambda}_{\b{z}} - f^{\lambda} &= (T_{\b{x}}+\lambda I)^{-1}g_{\b{z}} - f^{\lambda}
						  = (T_{\b{x}}+\lambda I)^{-1}g_{\b{z}} - (T_{\b{x}}+\lambda I)^{-1} (T_{\b{x}}+\lambda I) f^{\lambda}\nonumber\\
						  &=  (T_{\b{x}}+\lambda I)^{-1}\big[ g_{\b{z}} - (T_{\b{x}}+\lambda I) f^{\lambda} \big]
						  =  (T_{\b{x}}+\lambda I)^{-1}\big( g_{\b{z}} - T_{\b{x}} f^{\lambda} - \lambda f^{\lambda} \big)\nonumber\\
						  &=  (T_{\b{x}}+\lambda I)^{-1}\big( g_{\b{z}} - T_{\b{x}} f^{\lambda} - S_Kf_{\rho} + Tf^{\lambda} \big)\nonumber\\
						  &= (T_{\b{x}}+\lambda I)^{-1}\left( g_{\b{z}}  - S_Kf_{\rho} \right) + (T_{\b{x}}+\lambda I)^{-1}(T-T_{\b{x}})f^{\lambda}\nonumber\\
						  &= (T_{\b{x}}+\lambda I)^{-1}\left( g_{\b{z}}  - S_Kf_{\rho} \right) + (T_{\b{x}}+\lambda I)^{-1}(T-T_{\b{x}}) (T+\lambda I)^{-1}S_Kf_{\rho}.\label{eq:fz-flambda}
	    \end{align}
	    Let us rewrite $(T+\lambda I)^{-1}$ by the $(A+UV)^{-1} = A^{-1} - A^{-1}U\left(I+VA^{-1}U\right)^{-1}VA^{-1}$
	    operator Woodbury formula \citep[][Theorem~2.1, page 724]{ding08spectrum}
	    \begin{align*}
		(T+\lambda I)^{-1} &= (\lambda I+S_KS_K^*)^{-1}
					= (\lambda^{-1} I) - (\lambda^{-1} I)  S_K \left[I+S_K^*(\lambda^{-1}I)S_K\right]^{-1}S_K^* (\lambda^{-1} I)\\
					&= (\lambda^{-1} I) - \lambda^{-1}  S_K (\lambda I+\tilde{T})^{-1}S_K^*.
	    \end{align*}
	    By the derived expression for $(T+\lambda I)^{-1}$, we get $(T+\lambda I)^{-1}S_Kf_{\rho} =  \lambda^{-1}S_Kf_{\rho} - \lambda^{-1}  S_K (\lambda I+\tilde{T})^{-1}S_K^*S_Kf_{\rho}
						    = \lambda^{-1}S_K\big[f_{\rho} - (\tilde{T}+\lambda I)^{-1}S_K^*S_Kf_{\rho} \big]$.
	    Plugging this result to Eq.~\eqref{eq:fz-flambda}, introducing the $g_{\rho} = S_Kf_{\rho}$ notation, using the triangle inequality we get
	    \begin{eqnarray*}
		\lefteqn{\big\|\sqrt{T}\big(f^{\lambda}_{\b{z}} - f^{\lambda}\big)\big\|_{\H}=}\nonumber\\
		    &&\hspace*{-0.65cm}= \left\| \sqrt{T}(T_{\b{x}}+\lambda I)^{-1}\left\{ \left( g_{\b{z}}  - S_Kf_{\rho} \right) + (T-T_{\b{x}})  \lambda^{-1}S_K\left[f_{\rho} - (\tilde{T}+\lambda I)^{-1}S_K^*S_Kf_{\rho} \right]\right\}\right\|_{\H}\\
		    &&\hspace*{-0.65cm} \le \big\|\sqrt{T}(T_{\b{x}}+\lambda I)^{-1}\big\|_{\L(\H)}  \Big( \left\|g_{\b{z}}  - g_{\rho}\right\|_{\H} + \vphantom{\left\| T-T_{\b{x}}\right\|_{\L(\H)} \lambda^{-1} \left\| S_K\left[f_{\rho} - (\tilde{T}+\lambda I)^{-1}S_K^*S_Kf_{\rho} \right] \right\|_{\H}}
		      \left\| T-T_{\b{x}}\right\|_{\L(\H)} \lambda^{-1} \big\| S_K\big[f_{\rho} - (\tilde{T}+\lambda I)^{-1}S_K^*S_Kf_{\rho} \big] \big\|_{\H} \Big).
	    \end{eqnarray*}

\subsubsection{Proof of the bound on $\left\|g_{\b{z}}  - g_{\rho}\right\|_{\H}$}\label{proof:gz-grho:bound}
		      As is known $g_{\b{z}} = \frac{1}{l}\sum_{i=1}^l K_{\mu_{x_i}}y_i$ [see Eq.~\eqref{eq:f_zlambda}] and $g_{\rho} = \int_X K_{\mu_x}f_{\rho}(\mu_x)\d \rho_X(\mu_x)$ \citep[][Eq.~(23), page 344]{caponnetto07optimal}. Let 
		      $\xi_i = K_{\mu_{x_i}}y_i\in \H$ ($i=1,\ldots,l$). In this case $\E[\xi_i] = g_{\rho}$, $g_{\rho} - g_{\b{z}} = \E[\xi_i] - \frac{1}{l}\sum_{i=1}^l \xi_i$, and 
			$\left\|\xi_i\right\|_{\H}^2 =  \big\|K_{\mu_{x_i}}y_i\big\|_{\H}^2 \le \left\|K_{\mu_{x_i}}\right\|_{\L(Y,H)}^2 \left\|y_i\right\|_{Y}^2 
			\le B_K C^2  \Rightarrow \left\|\xi_i\right\|_{\H}\le C \sqrt{B_K} \hspace{0.01cm} \Rightarrow \hspace{0.01cm}
			  \E\big[\left\|\xi_i\right\|_{\H}^2\big] \le C^2 B_K$
		      using the boundedness of kernel $K$ ($\left\|K_{\mu_{x_i}}\right\|_{\L(Y,H)}^2\le B_K$) and the boundedness of output $y$ ($\left\||y\right\|_Y\le C$).
		      Applying the Bernstein inequality (see Lemma~\ref{lemma:Bernstein} with $\K=\H$, $B=2C\sqrt{B_K}$, $\sigma=C\sqrt{B_K}$) one gets that for any $\eta \in (0,1)$
		      \begin{align*}
			\Pr\left(\left\|g_{\b{z}}  - g_{\rho}\right\|_{\H} \le 2 \left( \frac{2C\sqrt{B}_K}{l}+ \frac{C\sqrt{B_K}}{\sqrt{l}} \right)\log\left(\frac{2}{\eta}\right) \right) \ge 1-\eta.
		      \end{align*}

\subsubsection{Proof of the bound on $\left\| T-T_{\b{x}}\right\|_{\L(\H)}$} \label{proof:T-Tx:bound}
  Let $\xi_i=T_{\mu_{x_i}}\in \L_2(\H)$ ($i=1,\ldots,l$), then
		      $\E[\xi_i] = T$,
		      $T-T_{\b{x}} = T - \frac{1}{l}\sum_{i=1}^lT_{\mu_{x_i}}$, $\left\|\xi_i\right\|_{\L_2(\H)} = \| T_{\mu_{x_i}}\|_{\L_2(\H)} \le B_K$,
		      $\E\left[\left\|\xi_i\right\|_{\L_2(\H)}^2\right] \le B_K^2$.
		  Applying the $\left\| T-T_{\b{x}}\right\|_{\L(\H)}\le \left\| T-T_{\b{x}}\right\|_{\L_2(\H)}$ relation [see Eq.~\eqref{eq:opnorm<=HSnorm}] and the Bernstein inequality (see Lemma~\ref{lemma:Bernstein} with $\K=\L_2(\H)$, $B=2B_K$, $\sigma=B_K$), we obtain that for any 
		  $\eta\in(0,1)$
		  \begin{align*}
		      \Pr\left( \left\| T-T_{\b{x}}\right\|_{\L(\H)} \le 2 \left(\frac{2B_K}{l}+\frac{\sigma}{\sqrt{l}}\right)\log\left(\frac{2}{\eta}\right) \right) \ge 1-\eta.
		  \end{align*}

\subsubsection{Proof of the decomposition of $\big\| S_K \big( f_{\rho} - (\tilde{T}+\lambda I)^{-1}S_K^*S_Kf_{\rho} \big) \big\|_{\H}^2$} \label{proof:S_Klong}
		Since $\left\|S_Ka\right\|_{\H}^2  = \left<S_Ka,S_Ka\right>_{\H} = \left<S_K^*S_Ka,a\right>_{\rho} = \big<\tilde{T}a,a\big>_{\rho}$ $(\forall a\in L^2_{\rho_X})$
		      by the definition of the adjoint operator and $\tilde{T}=S_K^*S_K$ [see Eq.~\eqref{eq:def:Ttilde}], we can rewrite the target term as
		      \begin{eqnarray*}
			  \lefteqn{\left\| S_K\left[f_{\rho} - (\tilde{T}+\lambda I)^{-1}S_K^*S_Kf_{\rho} \right] \right\|_{\H}^2 =}\\
			  &&\hspace*{3cm} = \left<\tilde{T}\left[f_{\rho} - (\tilde{T}+\lambda I)^{-1}S_K^*S_Kf_{\rho} \right], f_{\rho} - (\tilde{T}+\lambda I)^{-1}S_K^*S_Kf_{\rho} \right>_{\rho}\\
			  &&\hspace*{3cm}\le  \left\| \tilde{T}\left[f_{\rho} - (\tilde{T}+\lambda I)^{-1}S_K^*S_Kf_{\rho} \right] \right\|_{\rho} \left\| f_{\rho} - (\tilde{T}+\lambda I)^{-1}S_K^*S_Kf_{\rho} \right\|_{\rho},
		      \end{eqnarray*}
			where the CBS (Cauchy-Bunyakovsky-Schwarz) inequality was applied. Since
		      \begin{align}
			(S_K^*S_K+\lambda I)S_K^* &= S_K^*(S_KS_K^*+\lambda I) \hspace*{1cm}
			S_K^*(S_KS_K^*+\lambda I)^{-1}  &= (S_K^*S_K+\lambda I)^{-1}S_K^* \nonumber\\
			S_K^*(S_KS_K^*+\lambda I)^{-1} S_K &= (S_K^*S_K+\lambda I)^{-1}S_K^*S_K\label{eq:T-trafo1}\\
			S_K^*(T+\lambda I)^{-1} S_K &= (\tilde{T}+\lambda I)^{-1}\tilde{T}\label{eq:T-trafo2}
		      \end{align}
		      using Eq.~\eqref{eq:T-trafo1} and the analytical expression for $f^{\lambda}$ [see Eq.~\eqref{eq:analytical expr:f-lambda}] we have
		      \begin{align}
			  (\tilde{T}+\lambda I)^{-1}\tilde{T} f_{\rho} &= (\tilde{T}+\lambda I)^{-1}S_K^*S_K f_{\rho} = \left(S_K^*S_K+\lambda I\right)^{-1}S_K^*S_K f_{\rho}\nonumber\\
			  & = S_K^*(S_KS_K^*+\lambda I)^{-1} S_K f_{\rho} = S_K^* f^{\lambda}\label{eq:SKf-lambda:conversion}
		      \end{align}
		      and $\big\| S_K\big[f_{\rho} - (\tilde{T}+\lambda I)^{-1}S_K^*S_Kf_{\rho} \big] \big\|_{\H}^2 \le \big\| \tilde{T}\big[f_{\rho} - (\tilde{T}+\lambda I)^{-1}S_K^*S_Kf_{\rho} \big] \big\|_{\rho}
			\left\|S_K^*f^{\lambda} - f_{\rho}\right\|_{\rho}$.
		      
\subsubsection{Proof of the bound on $\left\|S_K^*f^{\lambda} - f_{\rho}\right\|_{\rho}$} \label{proof:flambda-frho-with-h}
	    Let us apply (i) the $Af - f = Af -f - q' + q' = (A-I)(f-q') + Aq' -q'$
	    relation with $A=(\tilde{T}+\lambda I)^{-1}\tilde{T}$, $f=f_{\rho}$ and $q'=S_K^*q$, where $q\in \H$ is an arbitrary element from $\H$, (ii) Eq.~\eqref{eq:SKf-lambda:conversion}
	    and (iii) the triangle inequality to arrive at
	    \begin{align*}
		\big\|S_K^*f^{\lambda} - f_{\rho}\big\|_{\rho} &= \big\| (\tilde{T}+\lambda I)^{-1}\tilde{T}f_{\rho} - f_{\rho}\big\|_{\rho}\\
		&= \big\| \big[(\tilde{T}+\lambda I)^{-1}\tilde{T} - I \big](f_{\rho} - S_K^*q) + (\tilde{T}+\lambda I)^{-1}\tilde{T} S_K^*q - S_K^*q\big\|_{\rho}\\
		&\le \big\| \big[(\tilde{T}+\lambda I)^{-1}\tilde{T} - I \big](f_{\rho} - S_K^*q)\big\|_{\rho}  + \big\|(\tilde{T}+\lambda I)^{-1}\tilde{T} S_K^*q - S_K^*q\big\|_{\rho}.
	    \end{align*}
	    Below we give upper bounds on these two terms. 

	    First, notice that $\mu_x\in X \mapsto \left\|K(\mu_x,\mu_x)\right\|_{\L(Y)}\le B_K$. This boundedness 
	    with the strong continuity of $K_{(\cdot)}$ imply \citep[][Proposition~12]{carmeli06vector} that $\H\subseteq C(X,Y)$, i.e., $K$ is a 
	    Mercer kernel. Since $K_{\mu_x}$ is a Hilbert-Schmidt operator for all $\mu_x\in X$ [see Eq.~\eqref{eq:bounded kernel}], it is also a compact operator ($\forall \mu_x\in X$). The compactness of $K_{\mu_x}$-s with 
	    the bounded and Mercer property of $K$ guarantees the boundedness of $S_K^{*}$ and that $\tilde{T}$ is a \emph{compact}, positive, self-adjoint operator \citep[][Proposition~3]{carmeli10vector}.\vspace{1mm}
	    
	    \noindent\tb{Bound on $\| [(\tilde{T}+\lambda I)^{-1}\tilde{T} - I ](f_{\rho} - S_K^*q)\|_{\rho}$}:
			Since $\tilde{T}$ is a compact positive self-adjoint operator, by the spectral theorem \citep[][Theorem~4.27, page~127]{steinwart08support} there exist an 
			$(u_i)_{i\in I}$ countable  ONB in $\cl{\Im(\tilde{T})}$, and $a_1\ge a_2\ge \ldots>0$ such that
			$\tilde{T}f = \sum_{i\in I}a_i \left<f,u_i\right>_{\rho}u_i$ $(\forall f\in L^2_{\rho_X})$
			and let $(v_j)_{j\in J}$ ($J$ is also countable by the separability\footnote{$L^2_{\rho_X}=L^2(X,\left.\Bo(H)\right|_X,\rho_X;Y)$ is isomorphic to $L^2(X,\left.\Bo(H)\right|_X,\rho_X;\R)\otimes Y$, where $\otimes$ is the tensor product of 
			Hilbert spaces. 
			The separability follows from that of $Y$ and $L^2(X,\left.\Bo(H)\right|_X,\rho_X;\R)$; the latter holds \citep[Proposition 3.4.5]{cohn13measure} since $\left.\Bo(H)\right|_X$ is countably generated since $X\subseteq H$ is separable. \label{footnote:separability-of-L2}} of $L^2_{\rho_X}$) an ONB in $\Ker(\tilde{T}^*)=\Ker(\tilde{T})$; $L^2_{\rho_X}=\cl{\Im(\tilde{T})} \oplus \Ker(\tilde{T})$. Thus,
			\begin{align*}
			     \big\| \big[(\tilde{T}+\lambda I)^{-1}\tilde{T} - I \big](f_{\rho} - S_K^*q)\big\|_{\rho}^2 &=
			       \sum_{i\in I}\Big(\frac{a_i}{a_i+\lambda}-1\Big)^2 \left<f_{\rho}-S_K^*q,u_i\right>_{\rho}^2 + \sum_{j\in J} \left<f_{\rho}-S_K^*q,v_j\right>_{\rho}^2\\
			      &\le \sum_{i\in I} \left<f_{\rho}-S_K^*q,u_i\right>_{\rho}^2 + \sum_{j\in J} \left<f_{\rho}-S_K^*q,v_j\right>_{\rho}^2
			      = \left\|f_{\rho}-S_K^*q\right\|_{\rho}^2
			\end{align*}
			exploiting the Parseval's identity and that $\big(\frac{\lambda_i}{\lambda_i+\lambda}-1\big)^2\le 1$.\vspace{1mm}

		\noindent\tb{Bound on $\|(\tilde{T}+\lambda I)^{-1}\tilde{T} S_K^*q - S_K^*q\|_{\rho}$}: By using Eq.~\eqref{eq:T-trafo2}, Eq.~\eqref{eq:S_K^*:norm preserving},
			   and Lemma~\ref{lemma:lin trafo-ed:ridge} ($M=T=S_KS_K^*$, $\K=\H$, $f=q$, $a=\frac{1}{2}$), the target quantity can be bounded as
		      \begin{align*}
			   \big\|(\tilde{T}+\lambda I)^{-1}\tilde{T} S_K^*q - S_K^*q\big\|_{\rho} &= \big\| S_K^*(T+\lambda I)^{-1} S_K S_K^*q - S_K^*q\big\|_{\rho}\nonumber\\
			    &= \big\| \sqrt{T}\left[(T+\lambda I)^{-1} S_K S_K^*q - q\right]\big\|_{\H}\\
			    &= \big\| \sqrt{T}\left[(T+\lambda I)^{-1} Tq - q\right]\big\|_{\H}
			    \le \max\left(1,\left\|T\right\|_{\L(\H)}\right)\lambda^{\frac{1}{2}} \left\|q\right\|_{\H}.
		      \end{align*}

	    Making use of the two derived bounds, we get $\left\|S_K^*f^{\lambda} - f_{\rho}\right\|_{\rho} \le  \left\|f_{\rho}-S_K^*q\right\|_{\rho} + \max\big(1,\left\|T\right\|_{\L(\H)}\big)\lambda^{\frac{1}{2}} \left\|q\right\|_{\H}$.

\subsection{Supplementary Lemmas} \label{sec:supplement}
In this section, we list two lemmas used in the proofs.

\subsubsection{Bernstein's inequality \citep[Prop.~2, p.~345]{caponnetto07optimal}} \label{lemma:Bernstein} %caponnetto07optimal: Proposition~2, page 13
Let $\xi_i$ ($i=1,\ldots,l$) be i.i.d.\ realizations of a random variable on a $(\Omega,\mathscr{A},P)$ probability space with values in a separable Hilbert space $\K$. If there exist $B>0$, $\sigma>0$ 
constants such that
 $\left\|\xi(\omega)\right\|_{\K}  \le \frac{B}{2}$ \text{ a.s.},  $\E\left[\left\|\xi\right\|_{\K}^2\right]  \le \sigma^2$, 
then for all $l\ge 1$ and $\eta\in (0,1)$ we have $$\Pr\left(\Big\|\frac{1}{l}\sum_{i=1}^l \xi_i -  \E[\xi_1]  \Big\|_{\K} \le 2\left(\frac{B}{l} + \frac{\sigma}{\sqrt{l}}\right) \log\left(\frac{2}{\eta}\right) \right)\ge 1-\eta.$$

\subsubsection{Lemma on  bounded, self-adjoint compact operators; \citet[][Proposition~A.2, page 39]{sriperumbudur14density}} \label{lemma:lin trafo-ed:ridge}
      Let $M$ be a bounded, self-adjoint compact operator on a separable Hilbert space $\K$. Let $a\ge 0$, $\lambda>0$, and $s\ge 0$. Let 
      $f\in\K$ such that $f\in \Im\left(M^{s}\right)$. If $s+a>0$, then
      \begin{align*}
	  \left\|M^a\left[(M+\lambda I)^{-1}Mf -f \right]\right\|_{\K} & \le \max\left(1,\left\|M\right\|^{s+a-1}_{\L(\K)}\right)\lambda^{ \min(1,s+a)}\left\|M^{-s}f\right\|_{\K}.
      \end{align*}
    Note: specifically for $s=0$ we have $\Im\left(M^{s}\right)=\Im\left(I\right)=\K$, in other words, there is no additional range space constraint.

\section{Discussion of Our Assumptions}\label{sec:assumptions-discussion}
 We give a short insight into the consequences of our assumptions (detailed in Section~\ref{sec:assumptions}) and present some concrete examples.
    \begin{itemize}
     \item \textbf{Well-definedness of $\rho$:}\label{discussion:rho:well-defined} The boundedness and continuity of $k$ imply the measurability of $\mu: (\M^+_1(\X),\Bo(\tau_w))\rightarrow (H,\Bo(H))$.  Let $\tau$ denote the open sets on $H=H(k)$, $\left.\tau\right|_X=\{A\cap X: A \in \tau\}$ the subspace topology on $X$, 
	and $\left.\Bo(H)\right|_X=\{A\cap X: A\in \Bo(H)\}$ the subspace $\sigma$-algebra on $X$. By noting \citep[Corollary~5.2.13]{schwartz98analyse} that
	$\Bo\left(\left.\tau\right|_X\right) = \left.\Bo(H)\right|_X = \{A\in  \Bo(H): A\subseteq X\}\subseteq \Bo(H)$, the H-measurability of $\mu$ guarantees the measurability of 
	$\mu: (\M^+_1(\X),\Bo(\tau_w))\rightarrow (X,\left.\Bo(H)\right|_X)$, and hence the well-definedness of  $\rho$, the measure induced by $\M$ on $X\times Y$; for further details see \citep[Section~A.1.1]{szabo15twostage}.\footnote{Note that the referred proof
        also holds for separable Hilbert $Y$, and by the simplified reasoning above the original $X\in \Bo(H)$ condition could be avoided.} 
	\item \textbf{Separability of $X$:} separability of $\X$ and the continuity of $k$ implies the separability of $H=H(k)$ \citep[Lemma 4.33, page~130]{steinwart08support}. Also, since $X\subseteq H$, $X$ is separable.
	\item \textbf{Finiteness of $B_k$:} \label{discussion:B_k-finite}If $\X$ is compact, then the continuity of $k$ implies $B_k < \infty$.
	\item \textbf{Finiteness of $B_K$, compact metricness of $X$:} Let $\X$ be a compact metric space. In this case $\M^+_1(\X)$ is also compact metric \citep[Theorem~6.4, page~55]{parthasarathy67probability}. Hence if 
	    $\mu:(\M^+_1(\X),\tau_w)\rightarrow H(k)$ is continuous\footnote{For example, if $k$ is universal, then $\mu$ metrizes the weak 
	      topology $\tau_w$ \citep[Theorem~23, page~1552]{sriperumbudur10hilbert}, hence $\mu$ is continuous.} (not just measurable), then $X$ is compact metric and thus 
	      by the H{\"o}lder property of $K_{(\cdot)}$, it is continuous implying that  $B_K<\infty$.
	\item \textbf{$K$ properties:} It is known \citep[page~339-340]{caponnetto07optimal} that
	      \begin{align}
		  K(\mu_a,\mu_b) &= K_{\mu_a}^* K_{\mu_b}, \quad (\forall \mu_a,\mu_b\in X) \label{eq:K-phiK}\\
		\left\| K_{\mu_a}\right\|_{\L(Y,\H)} &=  \left\|K_{\mu_a}^*\right\|_{\L(\H,Y)} \le \sqrt{B_K}, \quad (\forall \mu_a\in X). \label{eq:phiK:bounded}
	      \end{align}
	      \tb{Remark:} In terms of Eq.~\eqref{eq:K-phiK}, the Eq.~\eqref{eq:bounded kernel} assumption means that the $\{K(\mu_a,\mu_a)\}_{\mu_a\in X}$ operators are trace class, specifically they are compact operators.
      \item  \textbf{Separability of $\H$:} \label{discussion-of-assumptions:sepH} The separability of $X$ and the continuity of $K$ imply the separability of $\H$. Indeed, since $\mu_a\mapsto K_{\mu_a}$ is H{\"o}lder continuous w.r.t.\ the Hilbert-Schmidt norm it is also continuous. 
	      As a result it is continuous w.r.t.\ the operator norm, and thus also w.r.t.\ the strong topology. Using this property with the finiteness of $B_K$ the separability of $\H$ follows
		\citep[Proposition 5.1, Corollary 5.2]{carmeli06vector}.
      \item Our assumptions imply \citet{caponnetto07optimal}'s conditions (not considering the $\P(b,c)$ prior requirement). Indeed
		\begin{enumerate}[labelindent=0cm,leftmargin=*,topsep=0cm,partopsep=0cm,parsep=0cm,itemsep=0cm]
		    \item $Y$ is a separable Hilbert space by assumption; the same property also holds for $\H$ as we have seen.
		    \item The measurability of $(\mu_x,\mu_t)\mapsto \left<K_{\mu_x}w,K_{\mu_t}v\right>_{\H}$ for $\forall w,v\in Y$ is guaranteed by the continuity of $K_{(\cdot)}$ w.r.t.\ the strong topology.
		    \item We have $\int_{X\times Y}\left\|y\right\|_Y^2 \d \rho_X(\mu_x,y) \le \int_{X\times Y}C^2 \d \rho_X(\mu_x,y) = C^2<\infty$ due to the boundedness of $y$, and hence $\exists \Sigma>0, \exists M>0$ such that for $\rho_X\text{-almost }\mu_x\in X$ 
			    \begin{align}
				%\int_{Y}e^{\frac{\left\|y-f_{\rho}(\mu_x)\right\|_Y}{M}}-\frac{\left\|y-f_{\rho}(\mu_x)\right\|_Y}{M}-1 \d  \rho(y|\mu_x) \le \frac{\Sigma^2}{2M^2}. \label{eq:bounded noise}
				\int_Y \left\|y-f_{\rho}(\mu_x)\right\|_Y^m\d \rho(y|\mu_x) \le \frac{m! \Sigma^2 M^{m-2}}{2}\quad (\forall m\ge 2). \label{eq:bounded noise}
			      \end{align}
			    Indeed, by \citep[Eq.~(33)]{caponnetto07optimal} the Bernstein condition \eqref{eq:bounded noise} holds if $\left\|y-f_{\rho}(\mu_x)\right\|_Y \le \frac{M}{2}$,  $\int_Y \left\|y-f_{\rho}(\mu_x)\right\|_Y^2 \d \rho(y|\mu_x) \le \Sigma^2$. 
			    In our case using the boundedness of $y$, the regression function is also bounded and the same holds for $\left\|y-f_{\rho}(\mu_x)\right\|_{Y}$ by the triangle inequality: 
			    $\left\|y-f_{\rho}(\mu_x)\right\|_Y \le C + \left\|f_{\rho}\right\|_{\H}\sqrt{B_K}$; thus, 
			    $M =2(C + \left\|f_{\rho}\right\|_{\H}\sqrt{B_K})$, $\Sigma = \frac{M}{2}$ is a suitable choice.
		    \item The Polishness of $X\times Y$ was used by \citet{caponnetto07optimal} to assure the existence of $\rho(y|\mu_a)$; we guaranteed this existence under somewhat milder conditions (see footnote~\ref{footnote:rho(|)}). 
		\end{enumerate}
    \end{itemize}
    \tb{Real-valued outputs}: We now consider the specific case of $Y=\R$, when the following simplifications and results hold. By noting that in this case $Tr(K^*_{\mu_a}K_{\mu_a})=K(\mu_a,\mu_a)$, Eq.~\eqref{eq:bounded kernel} simplifies to the boundedness of kernel $K$ in the traditional sense
			\begin{align}
			      K(\mu_a,\mu_a) \le B_K\quad (\forall \mu_a\in X). \label{eq:K:bounded}
			\end{align}
		      Eq.~\eqref{eq:K:Lip} reduces to the H{\"o}lder continuity of the canonical feature map
		      $\Psi_K(\mu_c) := K(\cdot,\mu_c):X\rightarrow \H$, in other words $\exists L>0$, $h\in (0,1]$ such that
		      $\left\|\Psi_K(\mu_a) - \Psi_K(\mu_b)\right\|_{\H} \le L \left\|\mu_a - \mu_b\right\|_H^h,\quad \forall (\mu_a,\mu_b)\in X\times X$.
		      In case of a linear kernel, $K(\mu_a,\mu_b)=\left<\mu_a,\mu_b\right>_H$, $(\mu_a,\mu_b\in X)$, the H{\"o}lder continuity of $\Psi_K$ holds with 
			  $L=1$, $h=1$, and $B_K=B_k$ is a suitable choice. Evaluating the kernel $K$ at the empirical embeddings 
			    $\mu_{\hat{x}_i} = \int_{\X}k(\cdot,u)\d \hat{x}_i(u) =  \frac{1}{N}\sum_{n=1}^N k(\cdot,x_{i,n})\in H$
			yields the standard set kernel 
		      \begin{align*}
			    K(\mu_{\hat{x}_i},\mu_{\hat{x}_j}) &=  \left<\mu_{\hat{x}_i},\mu_{\hat{x}_j}\right>_H = \left<\frac{1}{N}\sum_{n=1}^N k(\cdot,x_{i,n}),\frac{1}{N}\sum_{m=1}^N k(\cdot,x_{j,m})\right>_H
								= \frac{1}{N^2}\sum_{n,m=1}^N k(x_{i,n},x_{j,m}) 
		      \end{align*}
		      by the bilinearity of $\left<\cdot,\cdot\right>_H$ and the reproducing property of $k$.\vspace{1mm}

		    \noindent\tb{Remark}: One can define many nonlinear kernels (see Table~\ref{tab:K examples}) on mean embedded distributions.
		    These kernels are the natural extensions to distributions of the Gaussian \citep{christmann10universal}, exponential, Cauchy, generalized t-student and inverse multiquadric kernels. 
		    If $\X$ is a compact metric space and $\mu$ is continuous, then the $\Psi_K$ canonical feature maps, associated to $K$-s in Table~\ref{tab:K examples}, can be shown to satisfy our H{\"o}lder continuity requirement [Eq.~\eqref{eq:K:Lip}]; for details, see \citep[Section~A.1.5-A.1.6]{szabo15twostage}. 

\begin{table}
  \begin{center}
  \begin{tabular}{@{}c@{\hspace{0.05cm}}c@{\hspace{0.05cm}}c@{\hspace{0.05cm}}c@{\hspace{0.05cm}}c@{}}
    \toprule
      $K_G$ & $K_e$ & $K_C$ & $K_t$ & $K_i$\\\midrule
      $e^{-\frac{\left\|\mu_a-\mu_b\right\|_H^2}{2\theta^2}}$ & $e^{-\frac{\left\|\mu_a-\mu_b\right\|_H}{2\theta^2}}$ & $\left(1+\left\|\mu_a-\mu_b\right\|_H^2 / \theta^2\right)^{-1}$ & $\left(1+\left\|\mu_a-\mu_b\right\|_H^{\theta}\right)^{-1}$ & $\left(\left\|\mu_a-\mu_b\right\|_H^2+\theta^2\right)^{-\frac{1}{2}}$\\
      $h=1$ & $h=\frac{1}{2}$ & $h=1$  & $h=\frac{\theta}{2}$ ($\theta\le 2$) & $h=1$\\ \bottomrule
  \end{tabular}
  \caption{Nonlinear kernels on mean embedded distributions: $K=K(\mu_a,\mu_b)$; $\theta>0$. For the H{\"o}lder continuity of $\Psi_K$, we assume that $\X$ is a compact metric space and $\mu$ is continuous.}
  \label{tab:K examples}
  \end{center}
\end{table}

\acks{We would like to thank the anonymous reviewers for their highly valuable, constructive suggestions to improve the manuscript.
This work was supported by the Gatsby Charitable Foundation, NSF grant 1247658, and DOE grant DE-SC001114. A part of the work was carried out while Bharath K. Sriperumbudur was a research fellow in the Statistical Laboratory, Department of Pure Mathematics and Mathematical Statistics at the University of Cambridge, UK.}

{%\small
\bibliography{14-510}} %=>natbib

\newpage

\begin{center}
  \Large \tb{Supplement}
\end{center}

\section{Proofs of the $\left\| f_{\b{z}}^{\lambda} \right\|_{\H}^2$ Bound, Theorem~\ref{conseq:conv-rate} and Theorem~\ref{conseq:L2rate:Assumption2b}}
This section contains the derivations of the $\left\| f_{\b{z}}^{\lambda} \right\|_{\H}^2$ bound (used in Theorem~\ref{theo1}; see Section~\ref{proof:fzlambda-bound}), Theorem~\ref{conseq:conv-rate} (Section~\ref{proof:conseq:conv-rate}) and Theorem~\ref{conseq:L2rate:Assumption2b} (Section~\ref{proof:L2rate:Assumption2b}).

\subsection{Bound on $\left\| f_{\b{z}}^{\lambda} \right\|_{\H}^2$}\label{proof:fzlambda-bound}
 Below we derive the stated Eq.~\eqref{f:zlambda:improved-bound} bound for $\left\| f_{\b{z}}^{\lambda} \right\|_{\H}^2$; it is guaranteed to hold under the conditions of the bounds for $S_1$ and 
$S_2$ obtained in \citep[Eq.~(48), above Eq.~(46), Eq.~(43)]{caponnetto07optimal}; see $(\dag)_1$, $(\dag)_2$, $(\dag)_3$ below.

Applying the triangle inequality and the definition of $\B(\lambda)$ we get
\begin{align*}
  \left\| f_{\b{z}}^{\lambda} \right\|_{\H} &\le  \left\|f_{\b{z}}^{\lambda} - f^{\lambda}\right\|_{\H}  + \left\| f^{\lambda} - f_{\rho}\right\|_{\H} + \left\| f_{\rho}\right\|_{\H} =  \left\|f_{\b{z}}^{\lambda} - f^{\lambda}\right\|_{\H}  + \sqrt{\B(\lambda)} + \left\| f_{\rho}\right\|_{\H}.
\end{align*}

$f_{\b{z}}^{\lambda} - f^{\lambda}$ can be decomposed \citep[page~347]{caponnetto07optimal} as 
\begin{align*}
  f_{\b{z}}^{\lambda} - f^{\lambda} &= (T_{\b{x}} + \lambda I)^{-1}(g_{\b{z}}-T_{\b{x}}f_{\rho}) + (T_{\b{x}} + \lambda I)^{-1} (T-T_{\b{x}})(f^{\lambda}-f_{\rho}).
\end{align*}
Thus
\begin{align*}
  \left\| f_{\b{z}}^{\lambda} - f^{\lambda}\right\|_{\H} &\le \underbrace{\left\| (T_{\b{x}} + \lambda I)^{-1}(g_{\b{z}}-T_{\b{x}}f_{\rho})\right\|_{\H}}_{=:(*)_1} + \underbrace{\left\|(T_{\b{x}} + \lambda I)^{-1} (T-T_{\b{x}})(f^{\lambda}-f_{\rho})\right\|_{\H}}_{=:(*)_2}.
\end{align*}
\begin{enumerate}[labelindent=0cm,leftmargin=*,topsep=0cm,partopsep=0cm,parsep=0cm,itemsep=0cm]
    \item 
	  The first term can be estimated as 
	  \begin{align*}
	  (*)_1 &\le  \left\| (T_{\b{x}} + \lambda I)^{-1}(T+\lambda I)^{\frac{1}{2}}\right\|_{\L(\H)} 
			\left\| (T+\lambda I)^{-\frac{1}{2}}(g_{\b{z}}-T_{\b{x}}f_{\rho})\right\|_{\H}\\
		&\stackrel{(\dag)_1}{\le}  \underbrace{\left\| (T_{\b{x}} + \lambda I)^{-1}(T+\lambda I)^{\frac{1}{2}}\right\|_{\L(\H)}}_{=:(*)_3} 2 \log\left(\frac{6}{\eta}\right)\left[\frac{1}{l}\sqrt{\frac{M^2 B_K}{\lambda}} + \sqrt{\frac{\Sigma^2\N(\lambda)}{l}} \right].
	  \end{align*}
	  By \citep[page~350]{caponnetto07optimal}
	  \begin{align*}
	    (*)_3 & = \left\| (T + \lambda I)^{-\frac{1}{2}} \left[ I - (T+ \lambda I)^{-\frac{1}{2}}(T-T_{\b{x}})(T+\lambda I)^{-\frac{1}{2}} \right]^{-1} \right\|_{\L(\H)} \\
		  & \le \left\| (T + \lambda I)^{-\frac{1}{2}}\right\|_{\L(\H)} \left\| \left[ I - (T+ \lambda I)^{-\frac{1}{2}}(T-T_{\b{x}})(T+\lambda I)^{-\frac{1}{2}} \right]^{-1} \right\|_{\L(\H)} \stackrel{(\dag)_2}{\le} \frac{2}{\sqrt{\lambda}}.
	  \end{align*}
    \item The second expression can be bounded as 
	\begin{align*}
	    (*)_2 &\le \left\|(T_{\b{x}} + \lambda I)^{-1}\right\|_{\L(\H)} \left\| (T-T_{\b{x}})(f^{\lambda}-f_{\rho})\right\|_{\H}
		  \stackrel{(\dag)_3}{\le}  \frac{2}{\lambda}\log\left(\frac{6}{\eta}\right) \left[ \frac{2B_K\sqrt{\B(\lambda)}}{l} + \sqrt{\frac{B_K \A(\lambda)}{l}}\right].
	\end{align*}
\end{enumerate}
To sum up, we obtained that 
\begin{align*}
   \left\| f_{\b{z}}^{\lambda} \right\|_{\H} &\le \frac{2}{\sqrt{\lambda}} 2 \log\left(\frac{6}{\eta}\right)\left[\frac{1}{l}\sqrt{\frac{M^2 B_K}{\lambda}} + \sqrt{\frac{\Sigma^2\N(\lambda)}{l}} \right]\\
	  &\quad + \frac{2}{\lambda}\log\left(\frac{6}{\eta}\right) \left[ \frac{2B_K\sqrt{\B(\lambda)}}{l} + \sqrt{\frac{B_K \A(\lambda)}{l}}\right]  + \sqrt{\B(\lambda)} + \left\| f_{\rho}\right\|_{\H}, 
\end{align*}
and hence by Eq.~\eqref{eq:norm-eq}
\begin{align*}
   \left\| f_{\b{z}}^{\lambda} \right\|_{\H}^2 &\le 6\left(\frac{16}{\lambda} \log^2\left(\frac{6}{\eta}\right)\left[\frac{M^2 B_K}{l^2\lambda} + \frac{\Sigma^2\N(\lambda)}{l} \right]\right.\\
	  &\quad \hspace*{0.6cm} \left. + \frac{4}{\lambda^2}\log^2\left(\frac{6}{\eta}\right) \left[ \frac{4B_K^2\B(\lambda)}{l^2} + \frac{B_K \A(\lambda)}{l}\right]  + \B(\lambda) + \left\| f_{\rho}\right\|_{\H}^2\right). 
\end{align*}

\subsection{Proof of Theorem~\ref{conseq:conv-rate}}\label{proof:conseq:conv-rate}
In the following $\lambda$ is chosen to match the 'bias' ($\lambda^c$) and a 'variance' (other) term in  $r(l,\lambda)$ [see Eq.~\eqref{eq:f1bb}], guarantee that the matched terms dominate and the constraints in 
Eq.~\eqref{eq:f1bb} are also satisfied; according to our assumptions $1<b$ and $c\in(1,2]$.

\begin{align}
 r(l,\lambda)  &= \frac{1}{l^{2+a}\lambda^3} + \frac{1}{l^a\lambda} + \frac{1}{l^{a+1}\lambda^{2+\frac{1}{b}}} + \lambda^c + \frac{1}{l^2\lambda} + \frac{1}{l\lambda^{\frac{1}{b}}}\rightarrow 0, 
      \text{ s.t. } l \lambda^{\frac{b+1}{b}}\ge 1,\hspace{0.1cm} l^a\lambda^2\ge 1 \label{eq:f1bb}.
\end{align}

\noindent\tb{Cases in terms of $a$ choice}:
\begin{itemize}[labelindent=0cm,leftmargin=*,topsep=0cm,partopsep=0cm,parsep=0cm,itemsep=0cm]
  \item $a\in \left(0,\frac{b(c+1)}{bc+1}\right]$: Since $\frac{b(c+1)}{bc+1}<2\Leftrightarrow b(1-c)<2$ ($\Leftarrow$ $b>1$, $c>1$) one has $\frac{1}{l^2\lambda}\le \frac{1}{l^a\lambda}$, and 
    \begin{align}
      r(l,\lambda)  &= \frac{1}{l^{2+a}\lambda^3} + \frac{1}{l^a\lambda} + \frac{1}{l^{a+1}\lambda^{2+\frac{1}{b}}} + \lambda^c + \frac{1}{l\lambda^{\frac{1}{b}}}\rightarrow 0, 
		      \text{ s.t. } l \lambda^{\frac{b+1}{b}}\ge 1,\hspace{0.1cm} l^a\lambda^2\ge 1 \label{eq:f1:a1}.
    \end{align}
    \begin{itemize}[labelindent=0cm,leftmargin=*,topsep=0cm,partopsep=0cm,parsep=0cm,itemsep=0cm]
      \item $(\boxed{1}=)\boxed{4}$: [order of the 1st and 4th terms in Eq.~\eqref{eq:f1:a1} are the matched and the first term will be discarded]: $\frac{1}{l^{2+a}\lambda^3} = \lambda^c \Leftrightarrow \lambda = l^{-\frac{2+a}{3+c}}$, and 
	    \begin{align}
	      r(l) &= l^{-a + \frac{2+a}{3+c}} + l^{-(a+1)+\frac{2+a}{3+c}\left(2+\frac{1}{b}\right)} + \boxed{l^{-\frac{(2+a)c}{3+c}}} + l^{-1+\frac{2+a}{(3+c)b}}\rightarrow 0, \label{eq:f1:a2}\\
		    &\quad \text{ s.t. }l^{1-\frac{2+a}{3+c}\frac{b+1}{b}}\ge 1, l^{a-2\frac{2+a}{3+c}} \ge 1.  \nonumber
	    \end{align}
	    \begin{itemize}[labelindent=0cm,leftmargin=*,topsep=0cm,partopsep=0cm,parsep=0cm,itemsep=0cm]
	      \item $\ovalbox{3}\ge \ovalbox{2}$ [the 3rd term dominates the 2nd in Eq.~\eqref{eq:f1:a2}]:  $-\frac{(2+a)c}{3+c} \ge - (a+1)+\frac{2+a}{3+c}\left(2+\frac{1}{b}\right) \Leftrightarrow (a+1)(c+3) \ge (2+a)(c+2+\frac{1}{b})  \Leftrightarrow a \left(1-\frac{1}{b}\right) \ge c+1+\frac{2}{b} \Leftrightarrow a \ge \frac{c+1+\frac{2}{b}}{1-\frac{1}{b}}>2$, which contradicts to 
	      $a \le \frac{b(c+1)}{bc+1}<2$.
	    \end{itemize}
      \item $(\boxed{2}=)\boxed{4}$: $\frac{1}{l^a\lambda} = \lambda^c \Leftrightarrow \lambda = l^{-\frac{a}{c+1}}$, and 
	      \begin{align*}
		    r(l) &= l^{-(2+a)+\frac{3a}{c+1}} +  l^{-(a+1)+\frac{a}{c+1}\left(2+\frac{1}{b}\right)} + \boxed{l^{-\frac{ac}{c+1}}} + l^{-1+\frac{a}{b(c+1)}}\rightarrow 0,\\
			  &\quad  \text{ s.t. } l^{1-\frac{a(b+1)}{b(c+1)}}\ge 1,\hspace{0.1cm} l^{a-\frac{2a}{c+1}}\ge 1.
	      \end{align*}
	      \begin{itemize}[labelindent=0cm,leftmargin=*,topsep=0cm,partopsep=0cm,parsep=0cm,itemsep=0cm]
		    \item $\ovalbox{3} \ge \ovalbox{1}$: $-\frac{ac}{c+1} \ge -(2+a)+\frac{3a}{c+1} \Leftrightarrow (2+a)(c+1)\ge a(c+3) \Leftrightarrow c+1\ge a$; this property holds because $c+1>2\ge a$.
		    \item $\ovalbox{3} \ge \ovalbox{2}$: $-\frac{ac}{c+1} \ge -(a+1)+\frac{a}{c+1}\left(2+\frac{1}{b}\right) \Leftrightarrow (a+1)(c+1) \ge a \left(c+2+\frac{1}{b}\right) \Leftrightarrow c+1 \ge a \left(1+\frac{1}{b}\right)\Leftrightarrow \frac{c+1}{1+\frac{1}{b}} \ge a$; this holds since $\frac{c+1}{1+\frac{1}{b}}>\frac{b(c+1)}{bc+1} \ge a$.
		    \item $\ovalbox{3}\ge \ovalbox{4}$:  $-\frac{ac}{c+1} \ge -1+\frac{a}{b(c+1)} \Leftrightarrow b\ge a$. Since $\frac{b(c+1)}{bc+1}\le b \Leftrightarrow 1\le b$, the required property holds by our $a\le \frac{b(c+1)}{bc+1}$ assumption.
		    \item Constraint-1 in $r(l)$: It is sufficient that $1-\frac{a(b+1)}{b(c+1)}> 0 \Leftrightarrow a< \frac{b(c+1)}{b+1}$; this holds since $a\le \frac{b(c+1)}{bc+1}<\frac{b(c+1)}{b+1}$. 
		    \item Constraint-2 in $r(l)$: It is enough to have $a-\frac{2a}{c+1} > 0 \Leftrightarrow c+1> 2$, which holds since $c> 1$.
	      \end{itemize}
	    To sum up, in this case $\lambda = l^{-\frac{a}{c+1}}$ and the rate is $r(l)=l^{-\frac{ac}{c+1}}\rightarrow 0$.
      \item $(\boxed{3}=)\boxed{4}$: $\frac{1}{l^{a+1}\lambda^{2+\frac{1}{b}}} = \lambda^c \Leftrightarrow \lambda = l^{-\frac{a+1}{2+\frac{1}{b}+c}}$, and 
	    \begin{align*}
		r(l) &= l^{-(2+a)+\frac{3(a+1)}{2+\frac{1}{b}+c}} + l^{-a+\frac{a+1}{2+\frac{1}{b}+c}} + \boxed{l^{-\frac{c(a+1)}{2+\frac{1}{b}+c}}} + l^{-1+\frac{a+1}{2+\frac{1}{b}+c}\frac{1}{b}}\rightarrow 0,\\
		&\quad \text{ s.t. }l^{1-\frac{a+1}{2+\frac{1}{b}+c}\frac{b+1}{b}}\ge 1, l^{a-\frac{2(a+1)}{2+\frac{1}{b}+c}} \ge 1.
	    \end{align*}
	    \begin{itemize}[labelindent=0cm,leftmargin=*,topsep=0cm,partopsep=0cm,parsep=0cm,itemsep=0cm]
	      \item $\ovalbox{3}\ge \ovalbox{2}$: $-\frac{c(a+1)}{2+\frac{1}{b}+c} \ge -a+\frac{a+1}{2+\frac{1}{b}+c} \Leftrightarrow a \left(2+\frac{1}{b}+c\right) \ge (a+1)(c+1) \Leftrightarrow a(1+\frac{1}{b}) \ge c+1 \Leftrightarrow a \ge \frac{c+1}{1+\frac{1}{b}}$.
		    However, $\frac{c+1}{1+\frac{1}{b}} > \frac{b(c+1)}{bc+1}$ ($\Leftrightarrow c>1$) which contradicts to our
		    $a\le \frac{b(c+1)}{bc+1}$ assumption.
	    \end{itemize}
      \item $(\boxed{5}=)\boxed{4}$: $\lambda^c = \frac{1}{l\lambda^{\frac{1}{b}}} \Leftrightarrow \lambda = l^{-\frac{1}{c+\frac{1}{b}}}$, and 
	  \begin{align*}
	      r(l) &= l^{-(2+a)+\frac{3}{c+\frac{1}{b}}} + l^{-a + \frac{1}{c+\frac{1}{b}}} + l^{-(a+1)+\frac{1}{c+\frac{1}{b}}\left(2+\frac{1}{b}\right)} + \boxed{l^{-\frac{c}{c+\frac{1}{b}}}}\rightarrow 0,\\
		  &\quad \text{ s.t. }l^{1-\frac{1}{c+\frac{1}{b}}\frac{b+1}{b}}\ge 1, l^{a -\frac{2}{c+\frac{1}{b}}}\ge 1.
	  \end{align*}
	  \begin{itemize}[labelindent=0cm,leftmargin=*,topsep=0cm,partopsep=0cm,parsep=0cm,itemsep=0cm]
	      \item $\ovalbox{4}\ge \ovalbox{2}$: $-\frac{c}{c+\frac{1}{b}} \ge -a + \frac{1}{c+\frac{1}{b}} \Leftrightarrow  a \ge \frac{c+1}{c+\frac{1}{b}}$. As we have seen ($\boxed{3}=\boxed{4}: \ovalbox{3}\ge \ovalbox{2}$) this contradicts to our $a$ choice.
	  \end{itemize}
    \end{itemize}
  \item $a\in \left(\frac{b(c+1)}{bc+1},\infty\right)$:
      \begin{align*}
      r(l,\lambda) &= \frac{1}{l^{2+a}\lambda^3} + \frac{1}{l^a\lambda} + \frac{1}{l^{a+1}\lambda^{2+\frac{1}{b}}} + \lambda^c + \frac{1}{l^2\lambda} + \frac{1}{l\lambda^{\frac{1}{b}}}\rightarrow 0, 
	    \text{ s.t. } l \lambda^{\frac{b+1}{b}}\ge 1,\hspace{0.1cm} l^a\lambda^2\ge 1.
      \end{align*}
      \begin{itemize}[labelindent=0cm,leftmargin=*,topsep=0cm,partopsep=0cm,parsep=0cm,itemsep=0cm]
	  \item $(\boxed{6} = )\boxed{4}$: Using the previous '$\boxed{5}=\boxed{4}$' case, $\lambda^c = \frac{1}{l\lambda^{\frac{1}{b}}} \Leftrightarrow \lambda = l^{-\frac{1}{c+\frac{1}{b}}}$, and 
		\begin{align*}
		    r(l) &= l^{-(2+a)+\frac{3}{c+\frac{1}{b}}} + l^{-a + \frac{1}{c+\frac{1}{b}}} + l^{-(a+1)+\frac{1}{c+\frac{1}{b}}\left(2+\frac{1}{b}\right)} + \boxed{l^{-\frac{c}{c+\frac{1}{b}}}} + l^{-2+\frac{1}{c+\frac{1}{b}}}\rightarrow 0,\\
			&\quad \text{ s.t. }l^{1-\frac{1}{c+\frac{1}{b}}\frac{b+1}{b}}\ge 1, l^{a -\frac{2}{c+\frac{1}{b}}}\ge 1.
		\end{align*}
		\begin{itemize}[labelindent=0cm,leftmargin=*,topsep=0cm,partopsep=0cm,parsep=0cm,itemsep=0cm]
		    \item $\ovalbox{4}\ge \ovalbox{1}$: $-\frac{c}{c+\frac{1}{b}} \ge -(2+a)+\frac{3}{c+\frac{1}{b}} \Leftrightarrow (2+a)\left(c+\frac{1}{b}\right)\ge c+3 \Leftrightarrow a\ge \frac{c+3}{c+\frac{1}{b}}-2 = \frac{b(c+3)}{bc+1}-2$. This requirement holds by our $a$ choice since 
		      $\frac{b(c+1)}{bc+1} \ge \frac{b(c+3)}{bc+1}-2\Leftrightarrow bc+1\ge b$ which is valid.
		    \item $\ovalbox{4}\ge \ovalbox{2}$: As we have seen (previous '$\boxed{5}=\boxed{4}$') this means $a\ge \frac{b(c+1)}{bc+1}$ which holds by our $a$ choice.
		    \item $\ovalbox{4}\ge \ovalbox{3}$: $-\frac{c}{c+\frac{1}{b}} \ge -(a+1)+\frac{1}{c+\frac{1}{b}}\left(2+\frac{1}{b}\right) \Leftrightarrow (a+1)\left(c+\frac{1}{b}\right) \ge 2+\frac{1}{b}+c \Leftrightarrow a \ge \frac{2+\frac{1}{b}+c}{c+\frac{1}{b}}-1 = \frac{2b+1+bc}{bc+1}-1$. This condition holds by our 
			  $a$ choice since $\frac{2b+1+bc}{bc+1}-1\le \frac{b(c+1)}{bc+1} \Leftrightarrow b\le bc$ which is valid.
		    \item $\ovalbox{4}\ge \ovalbox{5}$: $-\frac{c}{c+\frac{1}{b}} \ge -2+\frac{1}{c+\frac{1}{b}} \Leftrightarrow 2 \left(c+\frac{1}{b}\right) \ge c+1 \Leftrightarrow 2 \ge \frac{b(c+1)}{bc+1}$ what we have already established.
		    \item Constraint-1 in $r(l)$: It is enough to have $1-\frac{1}{c+\frac{1}{b}}\frac{b+1}{b}> 0 \Leftrightarrow 1 > \frac{b+1}{bc+1} \Leftrightarrow bc> b$, which holds.
		    \item Constraint-2 in $r(l)$: It is sufficient that $a -\frac{2}{c+\frac{1}{b}}> 0 \Leftrightarrow a > \frac{2}{c+\frac{1}{b}} = \frac{2b}{bc+1}$; this is satisfied because $a>\frac{b(c+1)}{bc+1}> \frac{2b}{bc+1}$.
		\end{itemize}
		To sum up, in this case $\lambda = l^{-\frac{1}{c+\frac{1}{b}}} = l^{-\frac{b}{bc+1}}$, the rate is $r(l)=l^{-\frac{bc}{bc+1}}$.	  
	  \item $(\boxed{1} = )\boxed{4}$: Using the previous '$\ovalbox{1} = \ovalbox{4}$' case,  $\frac{1}{l^{2+a}\lambda^3} = \lambda^c \Leftrightarrow \lambda = l^{-\frac{2+a}{3+c}}$, and 
	    \begin{align*}
	      r(l) &= l^{-a + \frac{2+a}{3+c}} + l^{-(a+1)+\frac{2+a}{3+c}\left(2+\frac{1}{b}\right)} + \boxed{l^{-\frac{(2+a)c}{3+c}}} + l^{-2 + \frac{2+a}{3+c}} + l^{-1+\frac{2+a}{(3+c)b}}\rightarrow 0,\\
		    &\quad \text{ s.t. }l^{1-\frac{2+a}{3+c}\frac{b+1}{b}}\ge 1, l^{a-2\frac{2+a}{3+c}} \ge 1.
	    \end{align*}
	      \begin{itemize}[labelindent=0cm,leftmargin=*,topsep=0cm,partopsep=0cm,parsep=0cm,itemsep=0cm]
		  \item $\ovalbox{3}\ge \ovalbox{5}$:    $-\frac{(2+a)c}{3+c} \ge -1+\frac{2+a}{(3+c)b} \Leftrightarrow (3+c)b\ge (2+a)(1+bc)\Leftrightarrow \frac{(3+c)b}{1+bc}-2\ge a$. However, 
		    $\frac{(3+c)b}{1+bc}-2\le \frac{b(c+1)}{bc+1} \Leftrightarrow b\le bc+1$ holds, thus the requirement can not be satisfied due to the $\frac{b(c+1)}{bc+1}<a$ assumption.
	      \end{itemize}
	  \item $(\boxed{2} = ) \boxed{4}$: Using the previous '$\boxed{2}=\boxed{4}$' case we have $\frac{1}{l^a\lambda} = \lambda^c \Leftrightarrow \lambda = l^{-\frac{a}{c+1}}$, and 
	      \begin{align*}
		    r(l) &= l^{-(2+a)+\frac{3a}{c+1}} +  l^{-(a+1)+\frac{a}{c+1}\left(2+\frac{1}{b}\right)} + \boxed{l^{-\frac{ac}{c+1}}} + l^{-2-\frac{a}{c+1}}+ l^{-1+\frac{a}{b(c+1)}}\rightarrow 0,\\
			  &\quad  \text{ s.t. } l^{1-\frac{a(b+1)}{b(c+1)}}\ge 1,\hspace{0.1cm} l^{a-\frac{2a}{c+1}}\ge 1.
	      \end{align*}
	      \begin{itemize}[labelindent=0cm,leftmargin=*,topsep=0cm,partopsep=0cm,parsep=0cm,itemsep=0cm]
		  \item $\ovalbox{3}\ge \ovalbox{5}$: $-\frac{ac}{c+1} \ge -1+\frac{a}{b(c+1)} \Leftrightarrow b(c+1)\ge a(1+bc) \Leftrightarrow \frac{b(c+1)}{bc+1}\ge a$; this contradicts to our $\frac{b(c+1)}{bc+1}< a$ choice.
	      \end{itemize}
	  \item $(\boxed{3} =) \boxed{4}$: Using the previous '$\boxed{3}=\boxed{4}$',  $\frac{1}{l^{a+1}\lambda^{2+\frac{1}{b}}} = \lambda^c \Leftrightarrow \lambda = l^{-\frac{a+1}{2+\frac{1}{b}+c}}$, and 
	    \begin{align*}
		r(l) &= l^{-(2+a)+\frac{3(a+1)}{2+\frac{1}{b}+c}} + l^{-a+\frac{a+1}{2+\frac{1}{b}+c}} + \boxed{l^{-\frac{c(a+1)}{2+\frac{1}{b}+c}}} + l^{-2+\frac{a+1}{2+\frac{1}{b}+c}} + l^{-1+\frac{a+1}{2+\frac{1}{b}+c}\frac{1}{b}}\rightarrow 0,\\
		&\quad \text{ s.t. }l^{1-\frac{a+1}{2+\frac{1}{b}+c}\frac{b+1}{b}}\ge 1, l^{a-\frac{2(a+1)}{2+\frac{1}{b}+c}} \ge 1.
	    \end{align*}
	    \begin{itemize}[labelindent=0cm,leftmargin=*,topsep=0cm,partopsep=0cm,parsep=0cm,itemsep=0cm]
		\item $\ovalbox{3}\ge \ovalbox{5}$: $-\frac{c(a+1)}{2+\frac{1}{b}+c} \ge -1+\frac{a+1}{2+\frac{1}{b}+c}\frac{1}{b} \Leftrightarrow 2b+1+bc \ge (a+1)(1+bc) \Leftrightarrow \frac{2b}{1+bc}-1\ge a$. However, $\frac{b(c+1)}{bc+1}\ge \frac{2b}{1+bc}-1\Leftrightarrow 2bc+1\ge b$ holds; thus, by the $a>\frac{b(c+1)}{bc+1}$  assumption the required property can not be satisfied.
	    \end{itemize}
	  \item $(\boxed{5} = )\boxed{4}$: $\lambda^c = \frac{1}{l^2\lambda} \Leftrightarrow \lambda = l^{-\frac{2}{c+1}}$, and 
	    \begin{align*}
		  r(l) &= l^{-(a+2)+\frac{6}{c+1}} + l^{-a+ \frac{2}{c+1}} + l^{-(a+1)+\frac{2}{c+1}\left(2+\frac{1}{b}\right)} + \boxed{l^{-\frac{2c}{c+1}}} +l^{-1+\frac{2}{b(c+1)}}\rightarrow 0,\\
		       &\quad \text{s.t. }l^{1-\frac{2(b+1)}{b(c+1)}}\ge 1,  l^{a-\frac{4}{c+1}}\ge 1.
 	    \end{align*}
	    \begin{itemize}[labelindent=0cm,leftmargin=*,topsep=0cm,partopsep=0cm,parsep=0cm,itemsep=0cm]
		\item $\ovalbox{4}\ge \ovalbox{5}$: $-\frac{2c}{c+1} \ge -1+\frac{2}{b(c+1)} \Leftrightarrow 1\ge \frac{2+2bc}{b(c+1)}\Leftrightarrow b\ge bc+2$ which does not hold.
	    \end{itemize}
      \end{itemize}
\end{itemize}

\subsection{Proof of Theorem~\ref{conseq:L2rate:Assumption2b}}\label{proof:L2rate:Assumption2b}
In the sequel we choose $\lambda$ by matching $2$ terms in $\sqrt{r(l,\lambda)}$ [Eq.~\eqref{eq:r:L2:rangespace-temp2}], guarantee that the matched terms dominate and the constraint in Eq.~\eqref{eq:r:L2:rangespace-temp2} holds; we proceed by matching the 'bias' ($\lambda^s$) and 'variance' (other) terms;
$s\in(0,1]$.
	\begin{align}
	    \sqrt{r(l,\lambda)} &= \frac{1}{l^a\lambda} +  \frac{1}{\lambda^{1-\frac{s}{2}}l^{\frac{1}{2}}} + \lambda^{s}\rightarrow 0,\text{ s.t. } l \lambda^2 \ge 1. \label{eq:r:L2:rangespace-temp2}
	\end{align}

\begin{itemize}[labelindent=0cm,leftmargin=*,topsep=0cm,partopsep=0cm,parsep=0cm,itemsep=0cm]
  \item $(\boxed{1} = )\boxed{3}$: $\frac{1}{l^a\lambda} = \lambda^{s} \Leftrightarrow \lambda = l^{-\frac{a}{s+1}}$, and
	\begin{align*}
	    \sqrt{r(l)} &= l^{-\frac{1}{2}+\frac{a}{s+1}(1-\frac{s}{2})} + \boxed{l^{-\frac{sa}{s+1}}} \rightarrow 0,\text{ s.t. } l^{1-\frac{2a}{s+1}} \ge 1.
	\end{align*}
	\begin{itemize}[labelindent=0cm,leftmargin=*,topsep=0cm,partopsep=0cm,parsep=0cm,itemsep=0cm]
	    \item $\ovalbox{2}\ge \ovalbox{1}$: $-\frac{sa}{s+1}\ge -\frac{1}{2}+\frac{a}{s+1}(1-\frac{s}{2}) \Leftrightarrow s+1\ge a(2+s) \Leftrightarrow \frac{s+1}{s+2} \ge a$.
	    \item Condition in $r(l)$: it is sufficient to have $1-\frac{2a}{s+1}> 0 \Leftrightarrow \frac{s+1}{2}> a$.
	\end{itemize}
	To sum up, if $a\le \frac{s+1}{s+2}\left(<\frac{s+1}{2}\right)$ then $\lambda = l^{-\frac{a}{s+1}}$ leads to the rate $\sqrt{r(l)}=l^{-\frac{sa}{s+1}}$; specifically, if $a=\frac{s+1}{s+2}$ then $\sqrt{r(l)}=l^{-\frac{s}{s+2}}$.
  \item $(\boxed{2} =) \boxed{3}$: $\frac{1}{\lambda^{1-\frac{s}{2}}l^{\frac{1}{2}}} = \lambda^{s} \Leftrightarrow \lambda = l^{\frac{-\frac{1}{2}}{1+\frac{s}{2}}=-\frac{1}{s+2}}$, and
	\begin{align*}
	    \sqrt{r(l)} &=  l^{-a  +\frac{1}{s+2}} + \boxed{l^{-\frac{s}{s+2}}} \rightarrow 0,\text{ s.t. } l^{1 -\frac{2}{s+2}} \ge 1.
	\end{align*}
	\begin{itemize}[labelindent=0cm,leftmargin=*,topsep=0cm,partopsep=0cm,parsep=0cm,itemsep=0cm]
	    \item $\ovalbox{2}\ge \ovalbox{1}$: $-\frac{s}{s+2} \ge -a  +\frac{1}{s+2} \Leftrightarrow a \ge \frac{s+1}{s+2}$.
	    \item Condition in $r(l)$: it is sufficient to have $1 -\frac{2}{s+2}> 0\Leftrightarrow s> 0$ which always holds.
	\end{itemize}
      To sum up,  if $\frac{s+1}{s+2} \le a$ then choosing $\lambda=l^{-\frac{1}{s+2}}$ the rate is $\sqrt{r(l)}=l^{-\frac{s}{s+2}}\rightarrow 0$.
\end{itemize}

\end{document}